\documentclass[a4paper,11pt,final]{article}

\usepackage{epsfig}
\usepackage[cmex10]{amsmath}
\usepackage{amssymb,amsthm}
\usepackage{mathtools}
\usepackage{nccmath}
\usepackage{soul}
\usepackage[normalem]{ulem}
\usepackage{cite}
\usepackage{todonotes}
\usepackage{tikz}
\usepackage{hyperref}
\usetikzlibrary{hobby,decorations.markings,arrows,intersections,shapes}
\usepackage{pgfplots}
\usepgfplotslibrary{fillbetween}
\pgfplotsset{compat=newest}
\usepgfplotslibrary{fillbetween}
\usepackage{verbatim}

\usepackage{lastpage,fancyhdr,graphicx}
\usepackage{epstopdf}
\usepackage{float}
\usepackage{tikz-cd}

\usepackage{amsthm,pifont}

\usepackage{caption}
\usepackage{subcaption}
\usepackage{multirow}

\newcommand{\SwSigj}{\sigma_{\! j}}
\newcommand{\SwSigF}{\sigma_{\! 1}}
\newcommand{\SwSigm}{\sigma_{\! m}}

\newcommand{\XSj}{{\scriptscriptstyle j}}

\begin{document}

% paper title: Must keep \ \\ \LARGE\bf in it to leave enough margin.
\title{\ \\ \LARGE\bf Model based MIN/MAX override control of centrifugal compressor
systems}

\author{ \Large Rico Schulze$^{1}$ and Hendrik Richter$^{2}$ \\ \\
	$^{1}$ AviComp Controls GmbH, Ostwaldstr. 4, \\ D–04329 Leipzig, Germany \\
	$^{2}$ HTWK Leipzig
University of Applied Sciences, Faculty of  Engineering,\\ Postfach 30 11 66, D–04251 Leipzig, Germany}

\maketitle

\begin{abstract} 

We consider an application-oriented  nonlinear control of centrifugal compressors. Industrial applications require the compressor system to adjust to variable process demands and to be restricted to the valid operation range (e.g. surge limit). We modify a compressor model of Gravdahl and Egeland to account for characteristic features of industrial compressors and combine the framework of nonlinear output regulation via the internal model principle with MIN/MAX-override control in order to implement trajectory tracking between given state constraints. Furthermore the switching scheme as well as the practical stability of the closed-loop MIMO system is analysed by the corresponding switched and impulsive error system. The override control is demonstrated by applying discharge pressure control, anti-surge control and maximum discharge pressure limitation.

\end{abstract}


\section{Introduction} % first section MUST be titled Introduction, and feature introductory text; do NOT change this title

Rotating compressors are an essential part of most plants of process industry. They are used to generate a pressure difference between the suction and the discharge side of the compressor, thus producing a continuous flow. Compressor operation is bound to a limited working range. The most important limitation is called surge limit and restricts operation to a certain minimum flow. Crossing the surge limit causes strong aerodynamic instabilities. This leads to a broad spectrum of negative consequences ranging from high vibrations to complete breakdown. Hence, most of the existing literature concerning control designs focuses on aerodynamic instabilities. Especially  so-called active control used for suppression of stall and surge by different design methods has been intensively studied. Among the methods,  Lyapunov-based designs (e.g. via backstepping) are most popular, e.g. \cite{Gravdahl:2002,Haddad:1999,Krstic:1998,Leonessa:2000a}. But also adaptive \cite{Blanchini:2002a} and high-gain control \cite{Blanchini:2002b} as well as bifurcation-based \cite{Liaw:1996} and optimum criteria-based methods \cite{Spakovszky:1999} have been used. Moreover, research was done on tracking fast setpoint changes for a coupled compressor/gas turbine system using flatness-based feedforward control \cite{Osmic:2014} and gain-scheduled decoupling control \cite{Schwung:2014}. An advantage of active control is the stabilization of unstable open-loop working points. This leads to an extension of the stable operation area in the closed-loop setup. However, an operation very close to the surge line is  explicitly allowed. This contradicts the safety and availability demands which industrial applications of plant operations  usually require. Therefore, surge avoidance rather than active suppression is typically applied in industrial applications. This corresponds to introducing state constraints in the controller design. Another important but rather rarely discussed aspect of industrial compressor operation is adjustment to variable process demands. Both aspects are taken into account by model predictive control (MPC). Linear MPC \cite{Bentaleb:2015,Budinis:2015,Cortinovis:2014,Cortinovis:2015} and nonlinear MPC \cite{Johansen:2002,Torrisi:2015,Vroemen:1999} have been applied to centrifugal compressor systems. 

Another way for meeting the requirements of industrial control of centrifugal compressors (tracking of constant as well as time-dependent trajectories together with introducing state constraints) is override control.~Override control is based on independent controllers that are organized in an override scheme and allows tracking of trajectories as long as no state constraint is violated. If any violation occurs, the control law is switched in order to return the state to the unrestricted domain. Despite the practical importance, there is little research effort on it. The most comprehensive work was done by Glattfelder and Schaufelberger \cite{Glattfelder:1988,Glattfelder:2003,Glattfelder:1983}. In particular, they investigated linear plant models in single-loop MIN/MAX override control. Another early work was conducted by Foss \cite{Foss:1981} on a linear gas turbine model. A somewhat different formulation of override control has been developed by Turner \& Postlethwaite \cite{Turner:2002}. Instead of a MIN/MAX structure switching between multiple sub-controllers, a single output violation compensator considering soft bounds for multivariable linear models has been designed. An extension of this approach to feedback linearizable systems has been reported by the same authors~\cite{Turner:2004}, while Herrmann~et~al.~\cite{Herrmann:2008} introduced strict bounds on linear multivariable plants. However, to the best of our knowledge only linear plant models have been considered for override control of centrifugal compressors so far. In this paper, a nonlinear model  is studied. 

A problem with override control of nonlinear compressor models is the need for a proof of stability of the switched system. To tackle the stability issue, we propose to combining MIN/MAX override control with the framework of internal model based nonlinear output regulation. This lead to a multiple controller formulation of MIN/MAX override control that provides practical stability as well as a simple setup and intuitive understanding in industrial environments.

The paper is structured as follows. In Section \ref{Sec:CompressorModeling} we give the centrifugal compressor model used for MIN/MAX override control, see also Appendix A for the full modelling procedure.~It expands the compressor model of Gravdahl and Egeland \cite{Gravdahl:1999b} to account for characteristic features of industrial compressors. The flow through the impeller and the diffuser is modelled separately so that the total effective passage length depends on the pressure ratio. In addition, the modified model  describes the jointed action of the adjustable positions of the guide vane and the blow-off valve while compensating  variable disturbances introduced by the process valve. This is a typical scenario in process control of industrial compressors. Section~\ref{Sec:OverrideControl} discusses the model based override control including a proof of practical stability. The application of the control scheme to the compressor model using nonlinear output regulation is demonstrated in Section~\ref{Sec:CompressorControl}. Section~\ref{Sec:Conclusions} gives conclusions and recommendations on further work.

\section{The centrifugal compressor model}\label{Sec:CompressorModeling}
Several models are available for describing stable and unstable compressor operations. Well-known examples for axial compressors are the models of Greitzer \cite{Greitzer:1976a} and Moore \& Greitzer \cite{Moore:1986b}. Based on Greitzer's lumped-parameter approach, turbocharger and centrifugal compressor models have been developed by Fink~et~al.~\cite{Fink:1992} and Gravdahl \& Egeland \cite{Gravdahl:1999b}, respectively. We consider a fixed-speed centrifugal compressor of industrial size. Since we focus on an application-oriented perspective we will expand the model of Gravdahl \& Egeland \cite{Gravdahl:1999b} to higher differential pressures that are common for centrifugal compressor stages of industrial size. This involves the flow through the impeller, the diffuser as well as the flow through control valves. Furthermore, we consider an application-oriented setup with a downstream process valve (PV) and a blow-off valve (BOV) that is used for the MIN/MAX anti-surge override control,  see Figure~\ref{figSys} of Appendix A for a schematic description.

To facilitate easy practical interpretability, the dynamics of the centrifugal compressor is described by the physical dimensional quantities, see the complete derivation of the model equations in Appendix A. In particular, we take the impeller velocity $c_2$ representing the flow though the compressor and the pressure ratio $\Pi=p_4/p_1$ accounting for the relation between plenum pressure and ambient pressure, as well as the guide vane position $r_{GV}$, the process valve position $r_{PV}$ and the blow-off valve position $r_{BOV}$. The compressor can be controlled by adjusting the positions of the guide vane $u_{GV}$ and the blow-off valve $u_{BOV}$. Adjusting the position of the process valve $u_{PV}$ is considered to model an external disturbance $z_{PV}:=u_{PV}$ to the compressor.

Using the Equations \eqref{eq:CentrifugalCompModel}, \eqref{eq:DynModelGV}, \eqref{eq:PlenumModelRe}, \eqref{eq:DynValveModels} and \eqref{eq:DynValveModels1} of Appendix A, and defining the state vector $x=(c_2,\Pi,r_{GV},r_{PV},r_{BOV})^T$, the control input vector $u=(u_{GV},u_{BOV})^T$ and the disturbance input $z=z_{PV}$, the compressor model becomes:
\begin{subequations} \label{eq:ComprModelFinal}
	\begin{align}
		\dot{x}_1 &= \frac{1}{L(x_2)} \left[Y_C(x_1,x_3) - k_1{\left(x_2^{r_k}-1\right)}\right] \label{eq:ComprModelFinala},\\
		\dot{x}_2 &= k_2 x_2^{r_k} \left[{x_1 - \left[k_{PV}Y_{PV}(x_4)+k_{BOV}Y_{BOV}(x_5)\right]\sqrt{{2k_1\left(x_2^{r_k}-1\right)}/r_k}}\right] \label{eq:ComprModelFinalb},\\
		\dot{x}_3 &= \frac{1}{\tau_{GV}}\left[u_1 - x_3\right] \label{eq:ComprModelFinalc},\\
		\dot{x}_4 &= \frac{1}{\tau_{PV}}\left[z - x_4\right] \label{eq:ComprModelFinald},\\
		\dot{x}_5 &= \frac{1}{\tau_{BOV}}\left[u_2 - x_5\right]\label{eq:ComprModelFinale},
	\end{align}
\end{subequations}
with the model parameters $k_1$, $k_2$, $r_k$ and the total effective passage length through the compressor $L(x_2)$:
\begin{subequations} \label{eq:ComprModelFinalParam}
	\begin{align}
	    k_1 &= R_S T_1, \quad k_2 = \kappa A_2/V, \quad r_k = (\kappa-1)/\kappa,\\	
        L(x_2) &= L_{12} + \frac{L_{23}}{x_2^{\frac{1}{\kappa}}\frac{A_3}{A_2}-1}\ln\left(x_2^{\frac{1}{\kappa}}\frac{A_3}{A_2}\right) + \frac{L_{34}}{x_2^{\frac{1}{\kappa}}\frac{A_3}{A_2}{\left(\frac{A_4}{A_3}-1\right)}}\ln\left(\frac{A_4}{A_3}\right).
	\end{align}
\end{subequations}
The nonlinear static compressor map $Y_C(x_1,x_3)$ is modeled by cubic polynomials according to a $\beta$-line interpolation along isolines, see Equation \eqref{eq:MapModelComplete} and Figure~\ref{fig:ComprCharacteristic} of Appendix A. 
The PV  characteristic $Y_{PV}(x_4)$, the BOV characteristic $Y_{BOV}(x_5)$ and the coefficients  $k_{PV}$ and $k_{BOV}$ are given by \eqref{eq:ValveFlowParameters}.

As we followed the basic modeling approach of Gravdahl \& Egeland, Equations~\eqref{eq:ComprModelFinala} and~\eqref{eq:ComprModelFinalb} are similar to \cite{Gravdahl:1999b} for constant rotational speed of the impeller. In contrast to Gravdahl \& Egeland, we have modeled the flow through the impeller and the diffuser separately. In combination with the isentropic valve flow model~(\ref{ThrottleModel}) the total effective passage length $L$ now depends on the pressure ratio $\Pi=x_2$. Hence, an increasing pressure ratio will lead to a reduction of the total effective passage length $L$ and vice versa. This accounts for the effect of changing dynamic responses of the system due to varying mass through the compressor. This effect is especially important for starting and stopping of the compressor as well as for compressor sections with a large operating range. In addition, we have taken into account the specific valve characteristics $Y_{PV}$ and $Y_{BOV}$. Hence, the model can be adopted to the application-specific designs of the valves where we use linear characteristic for the BOV and equal percentage for the PV  which complies with typical industrial setups.

Thus, from a modeling perspective, the model  \eqref{eq:ComprModelFinal} studied in this paper uses the physical variables flow and pressure ratio to describe the compressor operation, which is similarly to the Gravdahl \& Egeland model~\cite{Gravdahl:1999b},  but additionally accounts for the effect initiated by the jointed action of the adjustable positions of the guide vane and the blow-off valve. As the combination of guide vanes and blow-off valve is frequently used  in process control of industrial compressors, the model is more realistic from an application point of view.

\section{MIN/MAX override control} \label{Sec:OverrideControl}
\subsection{Problem statement}
\begin{figure}
\begin{center}
\includegraphics[width=10cm]{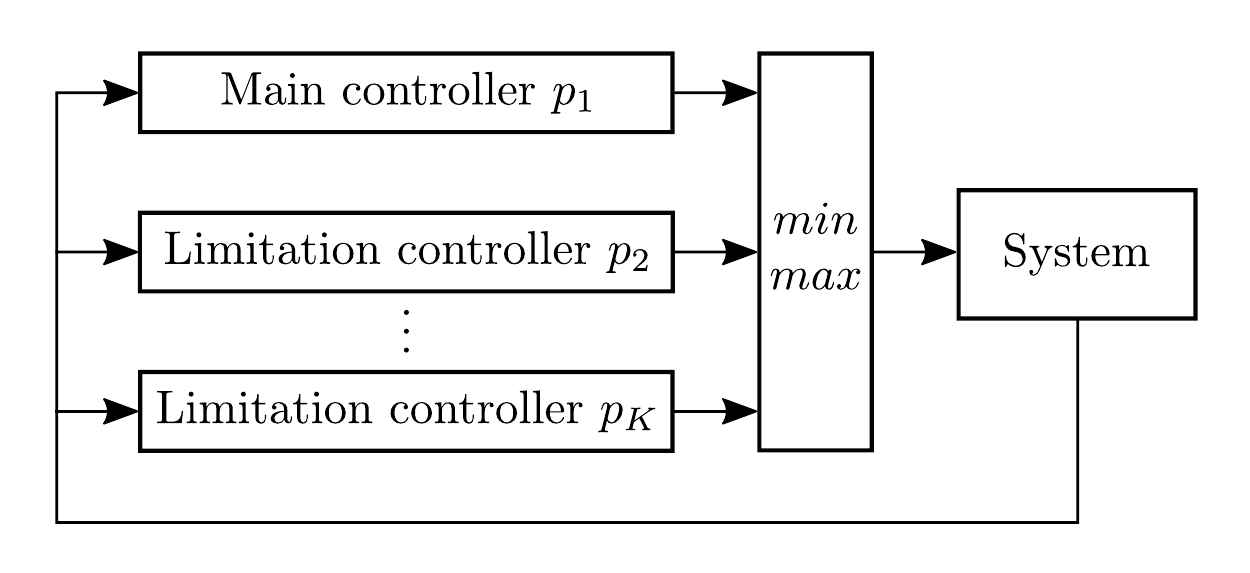}
\caption{General structure of MIN/MAX override consisting of a main controller and $K$ limitation controllers.}
\label{figOvrCntrl}
\end{center}
\end{figure}

Because active surge control is not applicable in an industrial environment we  adopt a surge avoidance strategy. The proposed structure of the MIN/MAX override control with state feedback is shown in Figure~\ref{figOvrCntrl}. From an application-oriented point of view this modular design is very advantageous as each controller can be designed and tuned independently. Roughly speaking such an controller design procedure corresponds to a divide-and-conquer approach. We consider one main controller and several limitation controllers. The main controller takes over the function of performance control, e.g. discharge pressure control. If an operating limit is reached the limitation controllers bound the main controller's output, e.g. for achieving anti-surge or maximum discharge pressure control. In total we have $K$ controllers which are indexed by $k$, i.e. $k \in \Omega\mathrel{\mathop:}=\{1,\ldots,K\}$. The compressor model~\eqref{eq:ComprModelFinal} is autonomous and input affine. At any time instant $t$ the compressor model~\eqref{eq:ComprModelFinal} is described by smooth vector fields $f,g_j:{\mathbb{R}}^{n}\rightarrow{\mathbb{R}}^{n}$: 
\begin{equation} \label{eq:AutoSys}
	\dot{x} = f(x) + \sum_{j=1}^{m}g_j(x)u_j,
\end{equation}
where $n$ is the dimension of the state space model, $m$ is the number of inputs and $u_j$ is the $j$-th input of the system. Now, we introduce state feedback controllers $p_k:{\mathbb{R}}^{n}\rightarrow{\mathbb{R}}^{m}$:
\begin{equation} \label{eq:ModeController}
	p_k(x) = \begin{bmatrix} p_{k1}(x) \\ \vdots \\ p_{km}(x) \end{bmatrix}.
\end{equation}
Here, $p_{kj}$ is the $j$-th sub-controller of the $k$-th controller. However, only one of the $K$ sub-controllers $p_{kj}$ can access the input $u_j$ at any time instant $t$ via the switching law $I_j$. In case of MIN/MAX override control each switching law $I_j:{\mathbb{R}}^K\rightarrow\Omega$ comprises of a combination of MIN/MAX selectors and generates a switching signal $\sigma_j \in \Omega$ via:
\begin{equation} \label{eq:SwitchSignal}
	\sigma_j(t)=I_j(p_{1j},\ldots,p_{Kj}).
\end{equation}
This implies that controller outputs are bounding each other. An override takes place if the controller output difference $\Delta p_{k_1k_2}^\XSj$ of two sub-controllers $p_{k_1j}(x)$ and $p_{k_2j}(x)$ ($k_1\neq k_2$) is zero, i.e.:
\begin{equation} \label{eq:RegulationDefOVR}
	\Delta p_{k_1k_2}^\XSj(t) = p_{k_1j}(t)-p_{k_2j}(t) = 0.
\end{equation}
An override from $p_{k_1j}(x)$ to $p_{k_2j}(x)$ will be marked by $p_{k_1j}\rightarrowtail p_{k_2j}$. We assume that the controllers are designed in a way to ensure $\Delta p_{k_1k_2}^\XSj(t)$ to be continuous. Otherwise, the state space cannot be divided into defined domains uniquely associated with a specific controller. This ensures that only one controller is active at a given time instant. When $\Delta p_{k_1k_2}^\XSj(t)$ is continuous the non-empty zero set $\mathcal{P}_{k_1k_2}^\XSj$ of~\eqref{eq:RegulationDefOVR} will be a continuous surface in the state space. The union of all zero sets $\mathcal{P}_{k_1k_2}^\XSj$ will be denoted as $\mathcal{P}$.

%The printed column width is 8.4 cm. Size the figures accordingly

As a control objective of the MIN/MAX override control is trajectory tracking, any controller $p_k$ must successfully achieve output regulation of system~\eqref{eq:AutoSys}. Therefore, the compressor model is expanded by the outputs (to be tracked) $h:{\mathbb{R}}^n\rightarrow{\mathbb{R}}^m\times{\mathbb{R}}^K$ with $h(x)=[h_1(x),\ldots,h_K(x)]$ and the setpoints $q:{\mathbb{R}}^n\rightarrow{\mathbb{R}}^m\times{\mathbb{R}}^K$ with $q(x)=[q_1(x),\ldots,q_K(x)]$ where $h_k(x)$ and $q_k(x)$ are given by:
\begin{equation} \label{eq:OutputSetpointDef}
	h_k(x) = \begin{bmatrix} h_{k1}(x) \\ \vdots \\ h_{km}(x) \end{bmatrix},\quad q_k(x) = \begin{bmatrix} q_{k1}(x) \\ \vdots \\ q_{km}(x) \end{bmatrix}.
\end{equation}
The choice of $h_k(x)$ and $q_k(x)$ depends on the control objective in the given application. Now, let us introduce the $k$-th controller tracking error $e_k:{\mathbb{R}}^n\rightarrow{\mathbb{R}}^{m}$:
\begin{equation} \label{eq:TrackingErrorDef}
		e_k = \begin{bmatrix} h_{k1}(x) - q_{k1}(x) \\ \vdots \\ h_{km}(x) - q_{km}(x) \end{bmatrix}.
\end{equation}
It is assumed that the tracking error $e_k$ is continuous. Hence, the non-empty zero set $\mathcal{E}_{k}$ of the tracking error~\eqref{eq:TrackingErrorDef} related to controller $p_{k}(x)$ will be a continuous surface in the state space. The union of all zero sets $\mathcal{E}_{k}$ will be denoted as $\mathcal{E}$. Using these definitions, system~\eqref{eq:AutoSys} can be rewritten:
\begin{subequations} \label{eq:AutoSysWithControl}
	\begin{align}
		\dot{x} &= f(x) + \sum_{j=1}^m g_j(x)\,p_{\SwSigj j}(x), \label{eq:AutoSysWithControla} \\
		e &= \begin{bmatrix} h_{\SwSigF 1}(x) - q_{\SwSigF 1}(x) \\ \vdots \\ h_{\SwSigm m}(x) - q_{\SwSigm m}(x) \end{bmatrix}, \label{eq:AutoSysWithControlb}
	\end{align}
\end{subequations}
where $e:{\mathbb{R}}^n\rightarrow{\mathbb{R}}^{m}$ is the tracking error of the active sub-controller. With respect to the above-named assumptions, the closed-loop system~\eqref{eq:AutoSysWithControl} is continuous despite the switching scheme. Thus, no sliding modes can occur which is a key advantage of the override control considered here. The steady state of system~\eqref{eq:AutoSysWithControl} is given by $e(t)=0$ and will be denoted as $\bar{x}(t)$; the steady state controller output is $p_k(x(t))\!\mid_{\,e(t)=0\;}=c_k(\bar{x}(t))$.

From an application point of view it is desirable to have a prescribed control dynamics. Therefore, we want the override control to satisfy the condition
\begin{equation} \label{cond:ControlDesignRegulation}
	\|e_{\SwSigj}(t)\|<\delta,\;t\in[t_{\delta,i},t_{i+1})
\end{equation}
after each switch, where $\delta\in\mathbb{R}_+$ is a sufficiently small positive constant and $t_{\delta,i}$ is a time instant with $t_i<t_{\delta,i}<t_{i+1}$. We are referring to condition~(\ref{cond:ControlDesignRegulation}) as the perfect regulation condition as is requires all elements of $e_{\SwSigj}(t)$ to be sufficiently small within a certain time interval after each switch. As we will see later this condition is sufficient for guaranteeing stability of the switched system as well as successful output regulation in a practical sense. Note that successful output regulation of any limitation controller also satisfies the given state constraints.

%defines a soft limit, i.e. overshoot is allowed during transient motion. However, during steady state all limits must be satisfied.

\subsection{Stability}\label{ssec:stability}
The stability of system~\eqref{eq:AutoSysWithControl} is determined by the stability of its corresponding error system. Hence, we shift an open-loop working point to the origin by separating the state and the input using the error coordinates $\tilde{x}, \tilde{u}$ and steady state coordinates $\bar{x}, \bar{u}$. Applying the separation $x(t) = \tilde{x}(t) + \bar{x}(t)$ and $u(t) = \tilde{u}(t) + \bar{u}(t)$ leads to the following switched impulsive system with state-depended switching:
\begin{subequations} \label{eq:AutoSysAsError}
	\begin{align} 
	\dot{\tilde{x}}(t) =\,& \tilde{f}_{\sigma}(\tilde{x}) + \sum_{j=1}^{m} \tilde{g}_{\SwSigj j}(\tilde{x})\bar{u}_{\SwSigj j} + \sum_{j=1}^{m} g_{\SwSigj j}(\tilde{x}) \tilde{u}_{{\SwSigj j}}, \quad \tilde{x}(t)\notin\mathcal{S}(t) \label{eq:AutoSysAsErrora},\\
		\tilde{x}(t_i^+) =\,& \tilde{x}(t_i^-) + \Delta \bar{x}(t_i), \quad \tilde{x}(t)\in\mathcal{S}(t) \label{eq:AutoSysAsErrorb},
	\end{align}
\end{subequations}
where $\Delta \bar{x}$ is the state reset and $\mathcal{S}(t)$ is the resetting set which describes the set of all states in the error state space where a switching (or override) takes place. Note, that $\mathcal{S}(t)$ is generated from the unified zero set $\mathcal{P}$. This zero set is now time-varying and exhibits impulsive changes during switchings. When $\tilde{x}(t)\in\mathcal{S}(t)$ then the trajectory $\tilde{x}(t)$ will jump from $\tilde{x}(t_i^-)$ to $\tilde{x}(t_i^+)$ at the time instant $t_i$ which is given by the resetting law~\eqref{eq:AutoSysAsErrorb}. It is assumed that $\tilde{x}(t_i^+)\notin\mathcal{S}(t_i^+)$. Hence, every switching instant $t_i$ is unique. Because of the continuous flow of system~\eqref{eq:AutoSysWithControl} in the error state space, we can evaluate the state reset as $\Delta \bar{x}(t_i) = \bar{x}(t^-_i) - \bar{x}(t^+_a)$, i.e. the state reset is given by the difference of steady states in its original coordinates. An inspection of the state reset reveals that in general $\exists t_i$ with $\|\Delta \bar{x}(t_i)\| > 0$. Using this, we make the following assumption and afterwards state the proposition.

\textbf{Assumption 1:} \textit{
Each sub-controller $p_k(x)$ with $k=1,\ldots,m$ guarantees asymptomatic stability of~\eqref{eq:AutoSysAsError} in an open neighborhood $\tilde{U}_A$ around $\tilde{x}=0$.}

This assumption appears reasonable as asymptotic stability of each sub-controller is a general requirement in industrial applications due to safety and performance reasons. It can be ensured by an appropriate design procedure, for instance as proposed by Huang~\cite{Huang:2004}, which  has also been used in controller design discussed in this paper, see Section 4.4.

\textbf{Proposition 1:}\textit{
Suppose Assumption 1 holds. If $\exists t_i$ with $\|\Delta \tilde{x}(t_i)\| > 0$ then there exists an open neighborhood $\tilde{U}_{\epsilon}$ around $\tilde{x}=0$ so that $\tilde{x}(t_i^-)\in \tilde{U}_{\epsilon} \Rightarrow \tilde{x}(t_i^+)\notin \tilde{U}_{\epsilon}$ follows $\forall\tilde{x}\in\tilde{U}_{\epsilon}$.}

\textbf{Proof:} \textit{Define an open neighborhood $\tilde{U}_{\epsilon}\mathrel{\mathop:} = \{\tilde{x}\in \mathbb{R}^n\mid \|\tilde{x}\| < \epsilon\}$ around $\tilde{x}=0$. From Assumption 1 it follows that there is some constant $\epsilon>0$ and a time instant $t_{\epsilon}$ for which the error state satisfies continuous motion $\|\tilde{x}(t)\|<\epsilon$ in the time interval $(t_{\epsilon},t_i)$. Now, choose $0<\epsilon\leq\|\Delta \tilde{x}(t^-_i)\|/2$. Then the error state $\tilde{x}(t)$ will leave $\tilde{U}_{\epsilon}$ during qualitative change at the time instant $t_i$ and the proposition follows.}$\qed$

Referring to the $\delta$-$\epsilon$ construction of Lyapunov stability, we can deduce an important consequence from Proposition 1. It follows that $\exists\epsilon>0$ so that $\forall\delta>0$, $\exists\|\tilde{x}(t_0)\|<\delta$ and $\exists t\geq t_0$ with $\|\tilde{x}(t)\|\geq\epsilon$. Hence, the switched impulsive error system~\eqref{eq:AutoSysAsError} cannot be stable in the sense of Lyapunov. As only certain $\varepsilon>0$ will satisfy the Lyapunov instability condition, we have to refer to explicit bounds resulting in the concept of practical stability. We will use the definition of~\cite{Lakshmikantham:1990}:

\textbf{Definition 1:} \textit{
System~\eqref{eq:AutoSysAsError} is said to be $(\lambda,A)$-practical stable if, given explicit bounds $(\lambda,A)$ with $0<\lambda<A$, then $\|\tilde{x}(t_0)\|<\lambda$ implies $\|\tilde{x}(t)\|<A$, $t\geq t_0$ for some $t_0\in\mathbb{R}_+$.}

We are interested in uniform practical stability where Definition 1 holds $\forall t_0\in\mathbb{R}_+$ \cite{Lakshmikantham:1990}. For considering practical stability of system~\eqref{eq:AutoSysAsError} another assumption on the switching events is required.

\textbf{Assumption 2:} \textit{
The sequence $\Sigma_{\Delta}=\{\|\Delta \bar{x}(t_1)\|,\|\Delta \bar{x}(t_2)\|,\|\Delta \bar{x}(t_3)\|,\ldots\}$ generated from the state resets at switching time instants $t_i$ is upper bounded by a sufficiently small constant $L\in\mathbb{R}_+$, i.e. $\|\Delta \bar{x}(t_i)\| < L$ $\forall i$.}

This assumption states that each state reset is bounded by an upper limit $L$. This limit must be set in a way that the region of stability for the given controller is not left after switching. We will assume that the controllers are designed in a way that this condition is met. Note that  Assumption 2 is obsolete if the controllers are designed to provide global stability.

Now, a theorem on the practical stability of the closed-looped switched and impulsive error system~\eqref{eq:AutoSysAsError} can be stated. The well-known dwell time stability approach as defined in \cite{Liberzon:2003} and \cite{Zhendong:2011} is applied.

\textbf{Theorem 1:} \textit{Suppose 
Assumption 1 and 2 hold. Consider the closed-loop error system~\eqref{eq:AutoSysAsError} with given initial conditions $\tilde{x}(t_0)\in \tilde{U}_{\lambda}$ and the open neighborhood $\tilde{U}_{\lambda}$ around $\tilde{x}=0$. If $t_{i+1}-t_i>\tau_i$ with sufficiently large dwell time $\tau_i$ holds for every override event, then error system~\eqref{eq:AutoSysAsError} is uniformly practically stable.}

\textbf{Proof:} \textit{To show practical stability according to Definition 1 two conditions must hold. The first condition is derived by examining the equation of impulsive motion~\eqref{eq:AutoSysAsErrorb}. If the error state satisfies $\tilde{x}(t_i^-)\in \tilde{U}_{\lambda}$ at any switching time instant $t_i$, then the state reset must be upper bounded in order to guarantee $\tilde{x}(t_i^+)\in \tilde{U}_A$. This is shown using the triangle inequality deduced from the equation of impulsive motion~\eqref{eq:AutoSysAsErrorb}:
\begin{equation}
	\|\tilde{x}(t_i^+)\|\leq\|\tilde{x}(t_i^-)\| + \|\Delta \bar{x}(t_i)\|<A.
\end{equation}
Now, if Assumption 2 holds with $\|\Delta \bar{x}(t_i)\|<L$, then we have $\|\tilde{x}(t_i^-)\| < \lambda$ with the upper limit $\lambda\in\mathbb{R}_+$ satisfying $\lambda < A-L$, i.e. Assumption 2 guarantees that the following condition holds for every switching event:
\begin{equation} \label{eq:FirstImplication}
	\|\tilde{x}(t_i^-)\| < \lambda \Rightarrow \|\tilde{x}(t_i^+)\| < A.
\end{equation}
The second condition is derived by examining the equation of continuous motion~\eqref{eq:AutoSysAsErrora}. If $\tilde{x}(t)\notin\mathcal{S}(t)$, then Assumption 1 guarantees that $\|\tilde{x}(t)\|$ will asymptotically decay to zero as time goes towards infinity. The $i$-th active controller must be activated for a minimum dwell time $\tau_i$ in order to reach $\tilde{U}_{\lambda}$ before the next switching event occurs. Hence, if switching is slow enough then $t_{i+1}-t_i>\tau_i$ holds and it follows that:
\begin{equation} \label{eq:SecondImplication}
	\|\tilde{x}(t_i^+)\| < A \Rightarrow \|\tilde{x}(t_{i+1}^-)\| < \lambda.
\end{equation}
Combining both conditions~\eqref{eq:FirstImplication} and~\eqref{eq:SecondImplication} leads to:
\begin{equation}
	\|\tilde{x}(t_1^-)\| < \lambda \Rightarrow \|\tilde{x}(t_1^+)\| < A \Rightarrow \|\tilde{x}(t_2^-)\| < \lambda \Rightarrow \|\tilde{x}(t_2^+)\| < A \Rightarrow \ldots
\end{equation}
which completes the proof for practical stability as required by  Definition 1.}$\qed$

Theorem 1 states that practical stability follows from (i.)~asymptotic stability of any controller within sufficiently large neighborhood of $\tilde{x}=0$,  (ii.)~bounded state resets and (iii.)~sufficiently slow switching.~Conversely, fast switching can cause instability which is a well-known phenomenon in industrial compressor override control. Note that Theorem 1 also satisfies the perfect regulation condition~\eqref{cond:ControlDesignRegulation} if $\lambda \leq \delta$, i.e. practical stability guarantees successful output regulation of the main controller and any limitation controller.

\section{Output regulation based on MIN/MAX override control} \label{Sec:CompressorControl}
\subsection{Internal model based output regulation}
We use the framework of the Byrnes-Isidori internal model principle for output regulation \cite{Isidori:1990} as a basis for designing the individual sub-controllers in the given MIN/MAX override scheme. It will be shown later that this framework can be applied to the compressor model \eqref{eq:ComprModelFinal} which can be described by system \eqref{eq:AutoSysWithControl} when MIN/MAX switching is considered. Following Isidori, system \eqref{eq:AutoSysWithControl} is expanded by the exo-state $v$:
\begin{subequations} \label{eq:AutoSysWithControlExo}
	\begin{align}
		\dot{x} &= f(x,v) + \sum_{j=1}^m g_j(x,v)\,p_{\SwSigj j}(x,v), \\
		\dot{v} &= s(v), \\
		e &= \begin{bmatrix} h_{\SwSigF 1}(x) - q_{\SwSigF 1}(v) \\ \vdots \\ h_{\SwSigm m}(x) - q_{\SwSigm m}(v) \end{bmatrix},
	\end{align}
\end{subequations}
where the exo-state $v=[w,d]^T$ consists of the controller setpoint exo-state $w$ and disturbance exo-state $d$. We assume that $s:{\mathbb{R}}^{o}\rightarrow{\mathbb{R}}^{o}$ is smooth and marginally stable as defined by \cite{Isidori:1995}. With the combined state $[x,v]^T$, the corresponding error system is again of the same structure as~\eqref{eq:AutoSysAsError}. Now, assume that $e = e_k$ holds, i.e. the $k$th sub-controller $p_k(x,v)$ accesses all control inputs of the system. An important feature of system~\eqref{eq:AutoSysWithControlExo} under $e = e_k$ is the existence of a center manifold with the mapping $x = \chi_k(v)$ for the $k$-th sub-controller. Following \cite{Isidori:1995}, this center manifold is locally invariant and guarantees exponential stability for the $k$-th sub-controller with $\|x(t)-\chi_k(v(t))\| \leq \beta_k \exp(-\alpha_k t) \|x(t_i)-\chi_k(v(t_i))\|$ ($\alpha_k,\beta_k\in\mathbb{R}_+$) in the time interval $(t_i,t_{i+1})$ if the initial tracking error satisfies $\|e_k(t_i^+)\|<A$ and the system is linear stable under $v(t)=0$, i.e. in the absence of the exo-state. Hence, the combined state $[x,v]^T$ converges to the steady state that is given by the mapping $x = \chi_k(v)$ and the objective of output regulation is reached for the $k$-th sub-controller $p_k(x,v)$. Hence, Assumption~1 will be satisfied as the statements above hold $\forall k\in \Omega$.

The approach just presented illustrates the main advantage of MIN/MAX override control: each sub-controller can be designed independently of each other. This is achieved by the specific MIN/MAX switching scheme together with the proof of stability from the last section. Hence, control design can be done $\forall k\in \Omega$ independent of the specific realization of the switching signal $\sigma$. It follows that the design for the switched system \eqref{eq:AutoSysWithControlExo} can be reduced to $K$ controller designs of the non-switched system:
\begin{subequations} \label{eq:AutoSysWithControlExok}
	\begin{align}
		\dot{x} &= f(x,w,d) + \sum_{j=1}^m g_j(x,w,d)\,p_{kj}(x,w,d), \label{eq:AutoSysWithControlExoka}\\
		\dot{w} &= s_w(w), \label{eq:AutoSysWithControlExokb}\\
		\dot{d} &= s_d(d), \label{eq:AutoSysWithControlExokc}\\
		e_k &= \begin{bmatrix} h_{k1}(x) - q_{k1}(w,d) \\ \vdots \\ h_{km}(x) - q_{km}(w,d) \end{bmatrix}\label{eq:AutoSysWithControlExokd},
	\end{align}
\end{subequations}
where $v=[w,d]^T$ has been applied as mentioned before. System \eqref{eq:AutoSysWithControlExok} is used for control design. We will refer to system \eqref{eq:AutoSysWithControlExok} as the combined compressor system. Before starting control design for the combined compressor system, some additional control requirements are discussed.

\subsection{Compressor control requirements} \label{subsec:ControlReq}

We consider discharge pressure control using an inlet or outlet guide vane as presented in Section \ref{Sec:CompressorModeling}. This is our main controller. The process valve (PV) is modeling the flow demand by the downstream process. Repositioning the PV will cause the compressor to change its working point. The main control objective is to keep the desired discharge pressure despite the disturbance introduced by the PV. 
Furthermore, compressor operation will be limited to a given maximum pressure as well as the given surge control line (SCL) using the BOV, where the BOV is opened to restrict compressor operation to the SCL during stationary operation. The BOV is the most effective way to avoid surge in industrial applications assuming that the handled gas can be vented to the environment. As long as the compressor is operating away from its limits, the BOV is closed. Hence, no gas is blown off at the discharge side and we have $k_{BOV}(x_5)=0$ and $u_{BOV}=0$. In this case the model is SISO. If the BOV is open the system is MIMO. It follows that in addition to the main control, two limiting controllers must be designed. The limiting controller will be referred to as maximum pressure limiter and anti-surge controller (ASC). The next step of control design is the trajectory generation for the setpoint and disturbance exo-systems \eqref{eq:AutoSysWithControlExokb} and \eqref{eq:AutoSysWithControlExokc}.

\subsection{Trajectory generation} \label{subsec:trajectory}
\begin{figure}
\begin{center}
\includegraphics[width=10cm]{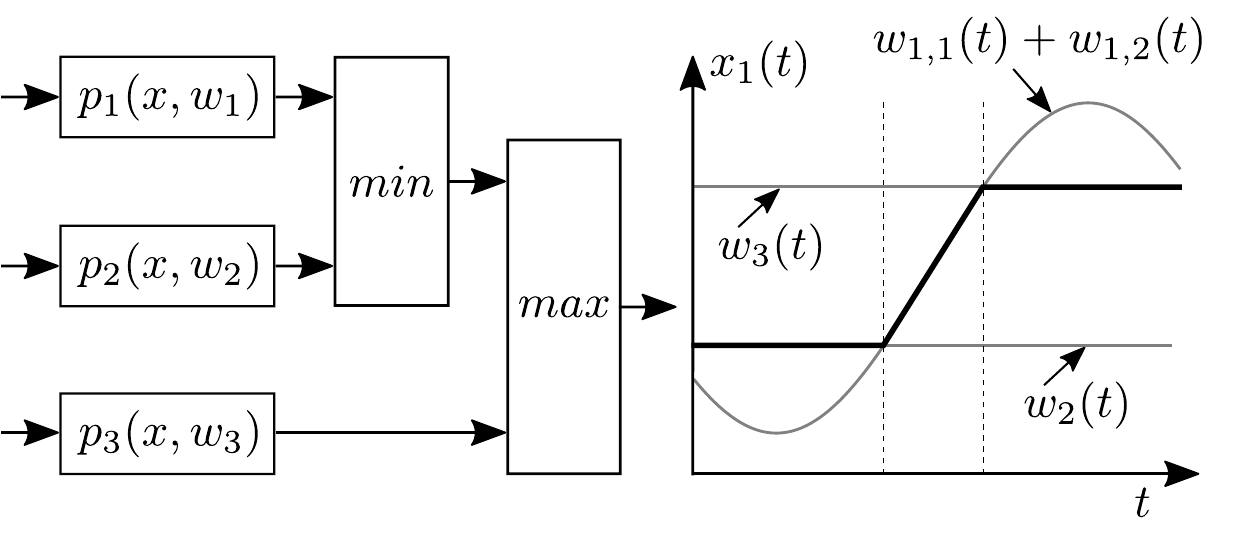}
\caption{Ramp generation by using MIN/MAX override of sinusoidal controller $p_1$ and two limitation controllers $p_2$ and $p_3$. An exemplary trajectory of $x_1(t)$ is plotted on the right hand side where the state $x_1(t)$ exhibits ramp-like motion (solid black line).}
\label{figRampGen}
\end{center}
\end{figure}

Industrial application of rotating compressor control mostly involve either constant setpoints or ramp-like setpoint motion. Both can be generated by using the concept of marginally stable exo-systems. However, the generation of a ramp signal is not possible for a single marginally stable exo-system \cite{Isidori:1995}. Therefore, we will approximate a ramp-like trajectory by using the suggested method of MIN/MAX override. Figure~\ref{figRampGen} illustrates the basic idea behind this approximation. A sinusoidal target trajectory, generated by the exo-system $\dot{w}_1 = \left[0, \omega w_{1,3}, -\omega w_{1,2}\right]^T$, moves between two constant limiting trajectories with $w_2$ and $w_3$ being their exo-states. Before reaching $t_i$ the system is in a steady state given by $w_2$ under control of $p_2(x,w_2)$. When reaching $t_i$ the MIN/MAX override structure selects $p_1(x,w_1)$ and starts to ramp. After reaching $t_{i+1}$ the MIN/MAX override structure  selects $p_3(x,w_3)$ and forces the system to approach the steady state given by $w_3$. To prevent the sinusoidal target trajectory from accessing the system again, $p_1(x,w_1)$ must be switched off once $p_3(x,w_3)$ is activated. Now, combining the ramp generator according to Figure~\ref{figRampGen}, we obtain the following exo-system:
\begin{subequations} \label{ComprCntrlExoSys}
	\begin{align}
		\dot{w}_1 &= \left[0, \omega w_{1,3}, -\omega w_{1,2}\right]^T, \\
		\dot{w}_2 &= 0, \\
		\dot{w}_3 &= 0,
	\end{align}
\end{subequations}
where the parameter $\omega$ determines the frequency of the sinusoidal main trajectory and at the same time defines the slope of the approximated ramp. Depending on the initialization of the corresponding exo-systems, the state either performs a positive or negative ramp-like motion. Note that exo-system~(\ref{ComprCntrlExoSys}) is marginally stable, which makes the presented approach a straightforward method to approximate ramp-like motion using marginally stable exo-systems.

%The printed column width is 8.4 cm. Size the figures accordingly

The method is now applied to the combined compressor model. The setpoint for the main controller (discharge pressure) and the ASC is provided by the above ramp generator with a sinusoidal trajectory that will be bounded by two constant trajectories with the exo-state variables $w_1$ to $w_5$. For the maximum pressure limiter we use a constant trajectory with just the exo-state variable $w_6$. This leads to the following setpoint exo-system $\dot{w} = s_w(w)$ with $s_w(w)$ given by:
\begin{equation} \label{ComprCntrlExoSysW}
    s_w(w) = \left[0, \omega_w w_{3}, -\omega_w w_{2}, 0,0,0\right]^T.
\end{equation}
Furthermore, we are modeling a sinusoidal motion for the disturbance that is introduced by the PV. Hence, the same structure as exo-system \eqref{ComprCntrlExoSys} is used for $\dot{d} = s_d(d)$ with $s_d(d)$ given by:
\begin{equation} \label{ComprCntrlExoSysD}
    s_d(d) = \left[0, \omega_d d_{3}, -\omega_d d_{2}\right]^T.
\end{equation}
This completes the definition of the exo-system. The next step of control design contains the compressor model analysis and the definition of the control law.

\subsection{Application to the compressor model} \label{subsec:appl}

\begin{figure}[!hb]
\begin{center}
	\includegraphics[width=\linewidth]{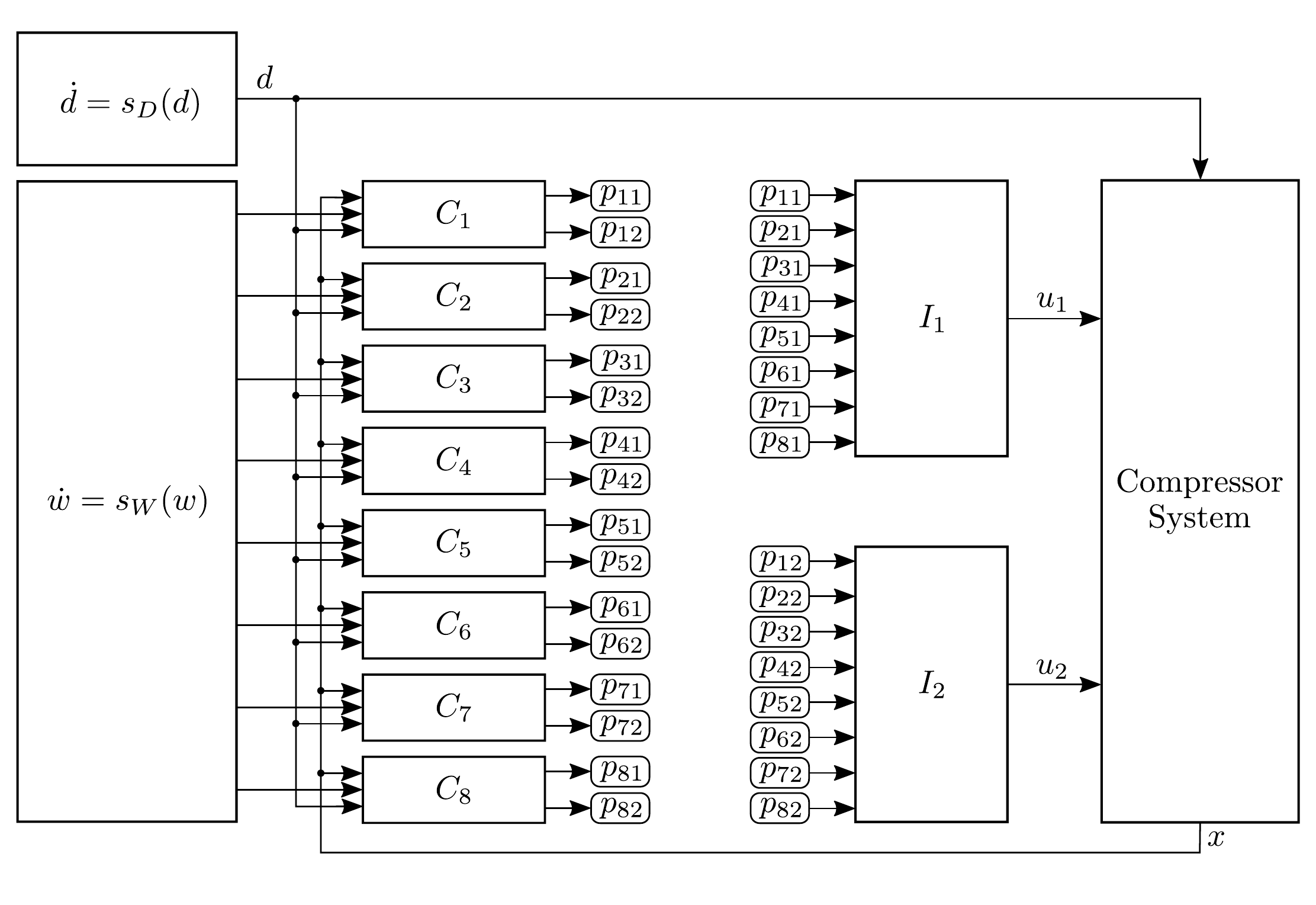}	\caption{Structure for internal model based output regulation controllers applied to compressor system model \eqref{eq:ComprModelFinal}.} \label{fig:ControlStructureCompressor}
\end{center}
\end{figure}

Following Section \ref{subsec:ControlReq}, several sub-controllers need to be designed for the SISO and MIMO domain. There is one main controller consisting of the sub-controllers $C_1$ (sine), $C_2$ (sine lower bound) and $C_3$ (sine upper bound). The ASC consists of the sub-controllers $C_5$ (sine), $C_6$ (sine lower bound) and $C_7$ (sine upper bound). The maximum pressure limitation controller consists of $C_4$ and $C_8$ both representing the maximum allowed pressure. Figure~\ref{fig:ControlStructureCompressor} illustrates  the control structure for the MIMO case.

%An examination of the total effective passage length $L$ reveals that its value only varies insignificantly compared to a suitably chosen averaged value $1/\tau_L$ as long as normal compressor operation is concerned. The maximum deviation does not exceed $3.2\%$ within the inspected range of $x_2=1.4\ldots2.2$ for the compressor under investigation. As no startup and shutdown scenarios are considered for control design and simulation, we will use the mean passage length $1/\tau_L$ in equation~\eqref{eq:ComprModelFinal}a. A likewise analysis of the factor $k_2 x_2^{r_k}$ in equation~\eqref{eq:ComprModelFinal}b shows that suitably chosen averaged value can be used as well. We will refer to the averaged value as $1/\tau_V$. The maximum deviation does not exceed $6.5\%$ within the inspected range of $x_2=1.4\ldots2.2$ for the compressor under investigation. 

Now, the control law is defined. We will use  $p_{kj}(x,w,d) = \bar{u}_{kj}(w,d) + G_{kj}(w,d)(x-\chi_k(w,d))$ for each sub-controller with the gain $G_{kj}(w,d)$ and the steady state setpoint $\bar{u}_{kj}(w,d)=p_{kj}(\chi_k(w,d),w,d)$. The state-feedback controller $p_{kj}(x,w,d)$ is a standard form that can solve the problem of output regulation and must be designed $\forall k\in\Omega$. Table~\ref{tab:MeasAndSPDefinition} contains the measurement equations and the setpoints used for control design, where $Y_{ASC}$ is the measurement equation for the ASC. We will assume a bivariate polynomial depending of the compressor flow $x_1$ and the pressure ratio $x_2$:
\begin{equation} \label{eq:ASCMeasEq}
    Y_{ASC}(x_1,x_2) = \sum_{n_1=0}^{N_1} \sum_{n_2=0}^{N_2} b_{n_1,n_2} x_1^{n_1} x_2^{n_2}
\end{equation}
with $N_1>0$. Furthermore, the following switching laws $I_1$ and $I_2$ according to Equation~\eqref{eq:SwitchSignal} are used:
\begin{subequations} \label{eq:SwitchingLaws}
	\begin{align}
		I_1 &= \max(\min(\max(p_{11},p_{21}),p_{31},p_{41}),\min(\max(p_{51},p_{61}),p_{71},p_{81})), \\
		I_2 &= \max(p_{12},p_{22},p_{32},p_{42},\min(\max(p_{52},p_{62}),p_{72},p_{82})).
	\end{align}
\end{subequations}
Before we can solve the output regulation problem for the proposed MIN/MAX override structure, it is verified whether the combined compressor system \eqref{eq:AutoSysWithControlExok} satisfies the necessary conditions for the control design presented by Huang \cite{Huang:2004}. Therefore, the relative vector degree and the zero dynamics are examined for the SISO and MIMO domain.

\renewcommand{\arraystretch}{1.1}
\setlength{\tabcolsep}{1.8mm}
\begin{table}[h!]
\centering
\begin{tabular}{cccc}
    \hline\hline
    Trajectory & SISO domain & MIMO domain & MIMO domain\\ 
    & & $j=1$ & $j=2$\\
    \hline
    \multirow{2}{*}{Sine}
    	& $h_1=x_2$			& $h_1=x_2$		    & $h_2=Y_{ASC}(x_1,x_2)$\\
    	& $q_{w,1}=w_1+w_2$	& $q_{w,1}=w_1+w_2$	& $q_{w,2}=Y_{ASC}(w_1,w_2)$\\
    \hline
    \multirow{2}{*}{Sine lower limit}
    	& $h_1=x_2$			& $h_1=x_2$		    & $h_2=Y_{ASC}(x_1,x_2)$\\
    	& $q_{w,1}=w_4$	& $q_{w,1}=w_4$	& $q_{w,2}=Y_{ASC}(w_4)$\\
    \hline
    \multirow{2}{*}{Sine upper limit}
    	& $h_1=x_2$			& $h_1=x_2$		    & $h_2=Y_{ASC}(x_1,x_2)$\\
    	& $q_{w,1}=w_5$	& $q_{w,1}=w_5$	& $q_{w,2}=Y_{ASC}(w_5)$\\
    \hline
    \multirow{2}{*}{Constant}
    	& $h_1=x_2$			& $h_1=x_2$		    & $h_2=Y_{ASC}(x_1,x_2)$\\
    	& $q_{w,1}=w_6$	& $q_{w,1}=w_6$	& $q_{w,2}=Y_{ASC}(w_6)$\\
    \hline\hline \\
\end{tabular}
\caption{Measurement equations and setpoint definition.}
\label{tab:MeasAndSPDefinition}
\end{table}
\renewcommand{\arraystretch}{1.5}

The relative vector degree for the combined compressor system \eqref{eq:AutoSysWithControlExok} is determined by using the combined notation:
\begin{subequations} \label{eq:B_ModelInStandardForm}
	\begin{align}
	\dot{x}_e =\;& f_e(x_e) + g_e(x_e) ,\\
	e =\;& h_e(x_e),
	\end{align}
\end{subequations}
with the combined state $x_e=[x,w,d]^T$ and the vector fields $f_e(x_e)=[f(x,d), s_W(w), s_D(d)]^T$, $g_e(x_e)=[g_1(x), 0_{m_W+m_D\times m_E}]^T$ and $h_e(x_e)=[h(x)-q_W(w)]$. For the SISO domain with $n=4$, $m_W=6$ and $m_D=3$ the relative degree is $r=3$. For the MIMO domain with $n=5$, $m_W=6$ and $m_D=3$. the relative vector degree is $r=\{2,2\}$. In both cases the combined internal dynamics is given by:
\begin{equation} \label{eq:B_CompleteZeroDynamics}
	\begin{bmatrix} \dot{x}_4 \\ \dot{w} \\ \dot{d} \end{bmatrix} = \begin{bmatrix} f_4(x_4,d) \\ s_W(w) \\ s_D(d) \end{bmatrix} = \begin{bmatrix} (q_D(d) - x_4)/\tau_{PV} \\ s_W(w) \\ s_D(d) \end{bmatrix}.
\end{equation}
The combined compressor system \eqref{eq:AutoSysWithControlExok} can be separated with $x^1=[x_1,x_2,x_3]^T$ and $x^2=[x_4]^T$ in the SISO domain as well as $x^1=[x_1,x_2,x_3,x_5]^T$ and $x^2=[x_4]^T$ in the MIMO domain. The combined internal dynamics \eqref{eq:B_CompleteZeroDynamics} does neither depend on $x^1$ nor on the inputs $u_1$ and $u_2$. Hence, the combined internal dynamics \eqref{eq:B_CompleteZeroDynamics} is also the combined zero dynamics of the combined compressor system. The output regulation problem can be solved if the subsystem specified by the $\dot{x}_4=f_4(x_4,d)$  of the combined zero dynamics in Equation \eqref{eq:B_CompleteZeroDynamics} is linearly stable \cite{Huang:2004}. As the subsystem specified by  $\dot{x}_4=f_4(x_4,d)$ is a first order low pass filter system, its characteristic polynomial is given by $\lambda \tau_{PV} + 1$. Assuming $\tau_{PV}>0$, exponential stability is given for $\dot{x}_4=f_4(x_4,d)$. Hence, the combined compressor system \eqref{eq:AutoSysWithControlExok} is minimum phase and we can apply the control design algorithm of Huang \cite{Huang:2004} to solve the center manifold $\chi_k(w,d)$ $\forall k\in \Omega$ and the corresponding steady state setpoint $\bar{u}_{kj}(w,d)=p_{kj}(\chi_k(w,d),w,d)$. Requirements on the control accuracy or on the tracking speed can be achieved by adjusting the control gains $G_{kj}(w,d)$. Each controller can be adjusted separately. As allready mentioned, this is considered an advantage of the presented modular design of a MIN/MAX override control.

\subsection{Simulation results} \label{subsec:sim}
\begin{figure}[!h]
    \begin{subfigure}{0.485\textwidth}
    	\includegraphics[width=\linewidth]{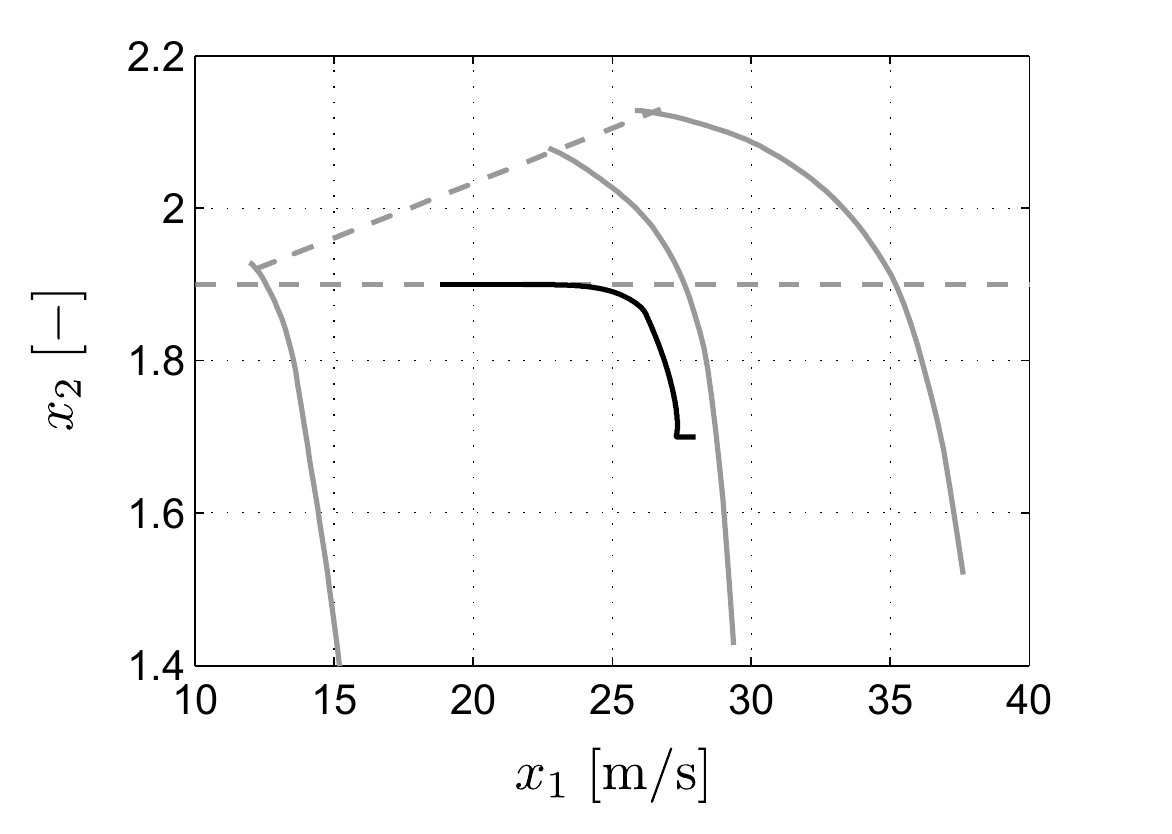}	\caption{Simulation results for scenario 1.}
        \label{fig:SimWPMotionS1}
    \end{subfigure}\hfill
    \begin{subfigure}{0.485\textwidth}
        \includegraphics[width=\linewidth]{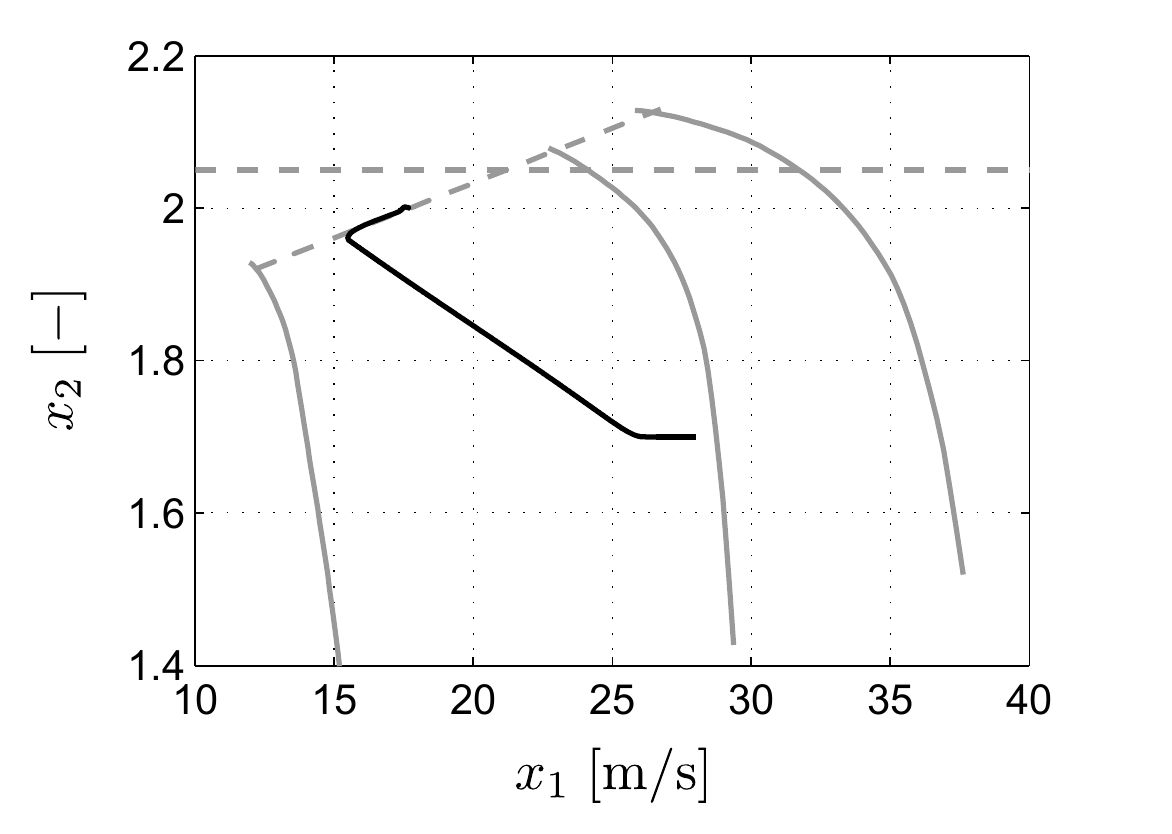}	\caption{Simulation results for scenario 2.} \label{fig:SimWPMotionS2}
    \end{subfigure}\hfill

	\caption{Compressor working point motion displayed in the $x_1$-$x_2$-plane (solid black) together with surge control line and maximum pressure limit (both dashed gray).}
	\label{fig:SimWPMotion}
\end{figure}

\begin{figure}[!hb]
    \begin{subfigure}[t]{1\textwidth}
        \begin{subfigure}{0.485\textwidth}
        	\includegraphics[width=\linewidth]{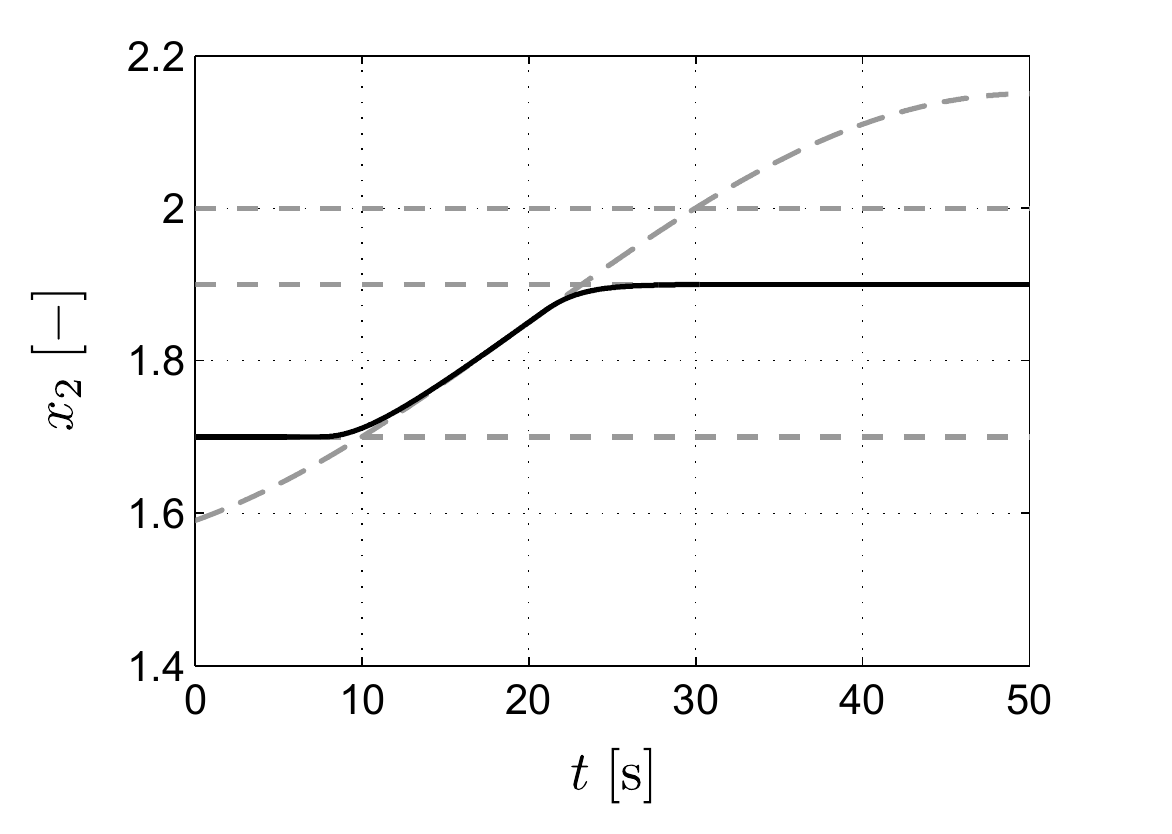}	\caption{Pressure ratio (solid black) together with defined trajectories (dashed gray) for scenario 1.}
            \label{fig:SimPiS1}
        \end{subfigure}\hfill
        \begin{subfigure}{0.485\textwidth}
            \includegraphics[width=\linewidth]{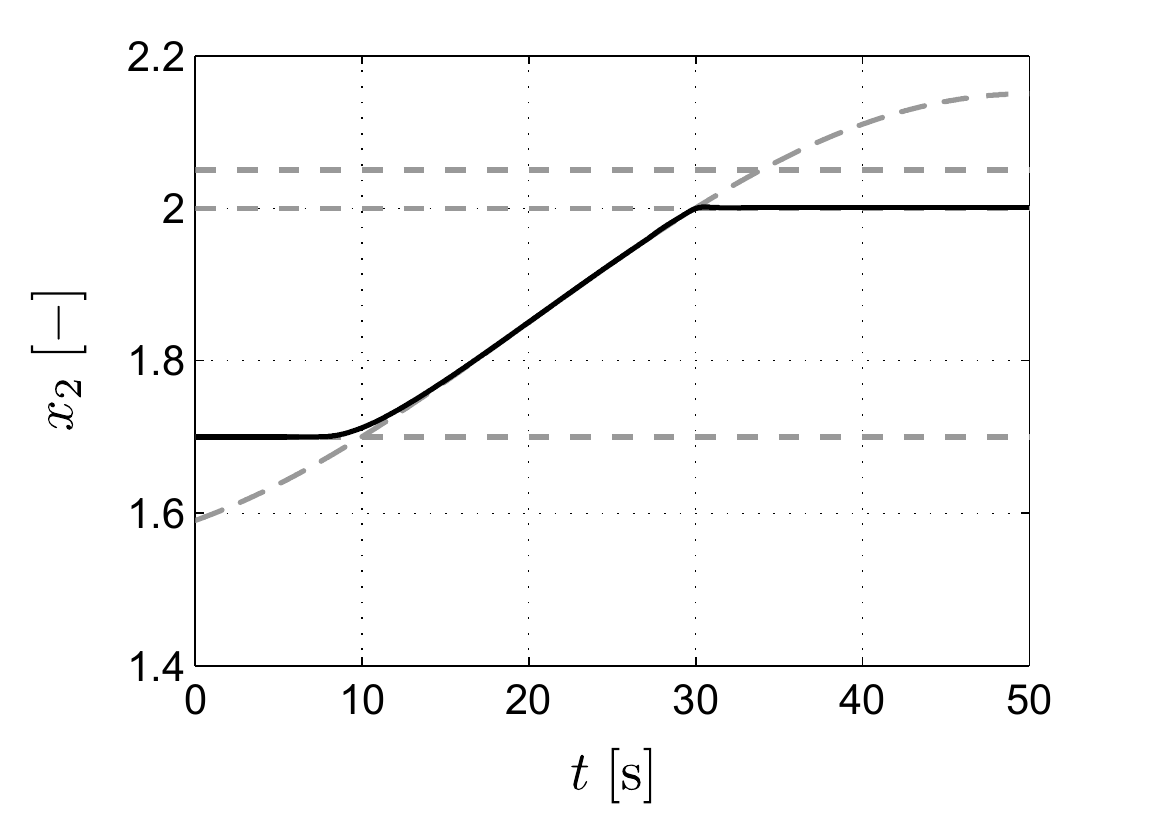}	\caption{Pressure ratio (solid black) together with defined trajectories (dashed gray) for scenario 2.} \label{fig:SimPiS2}
        \end{subfigure}\hfill
    \end{subfigure}\hfill

    \begin{subfigure}[t]{1\textwidth}
        \begin{subfigure}{0.485\textwidth}
        	\includegraphics[width=\linewidth]{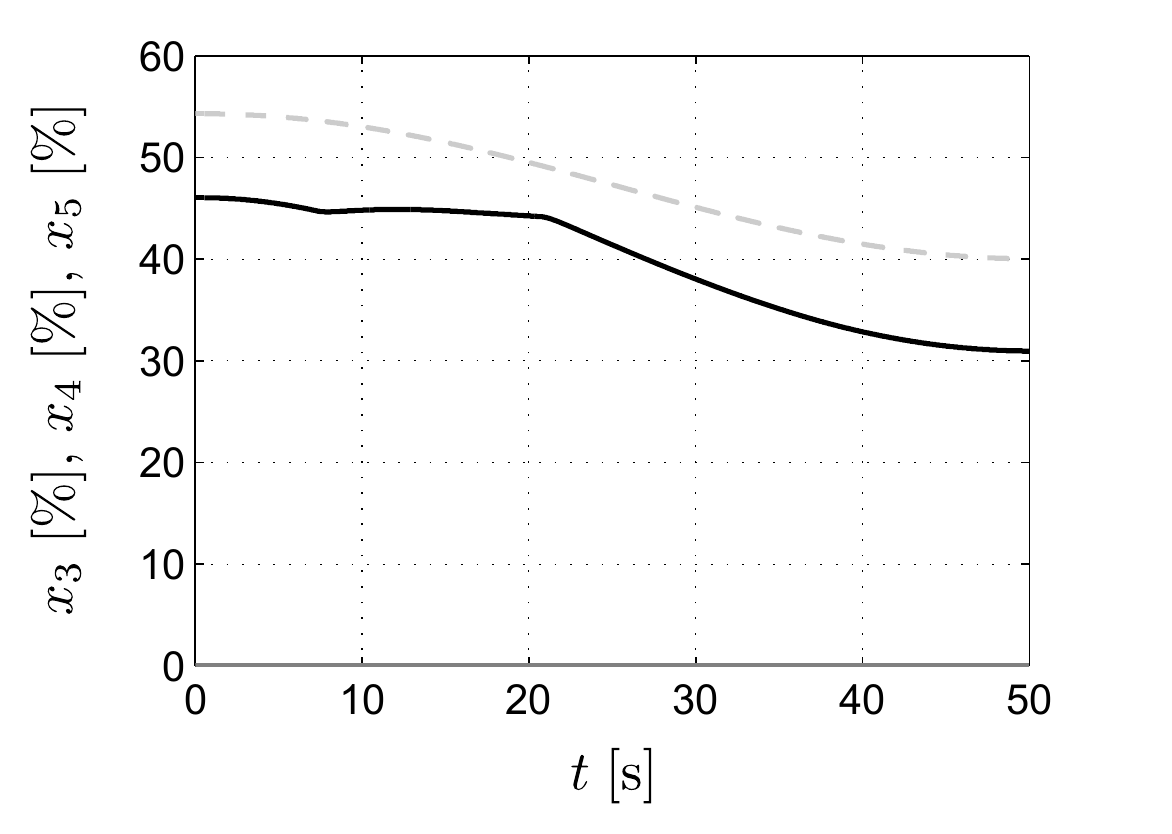}	\caption{Guide vane (GV) position $x_3$ (solid black), process valve (PV) position $x_4$ (dashed light grey) and blow-off valve (BOV) position $x_5$ (solid dark grey) for scenario 1.  }
            \label{fig:SimActuatorsS1}
        \end{subfigure}\hfill
        \begin{subfigure}{0.485\textwidth}
            \includegraphics[width=\linewidth]{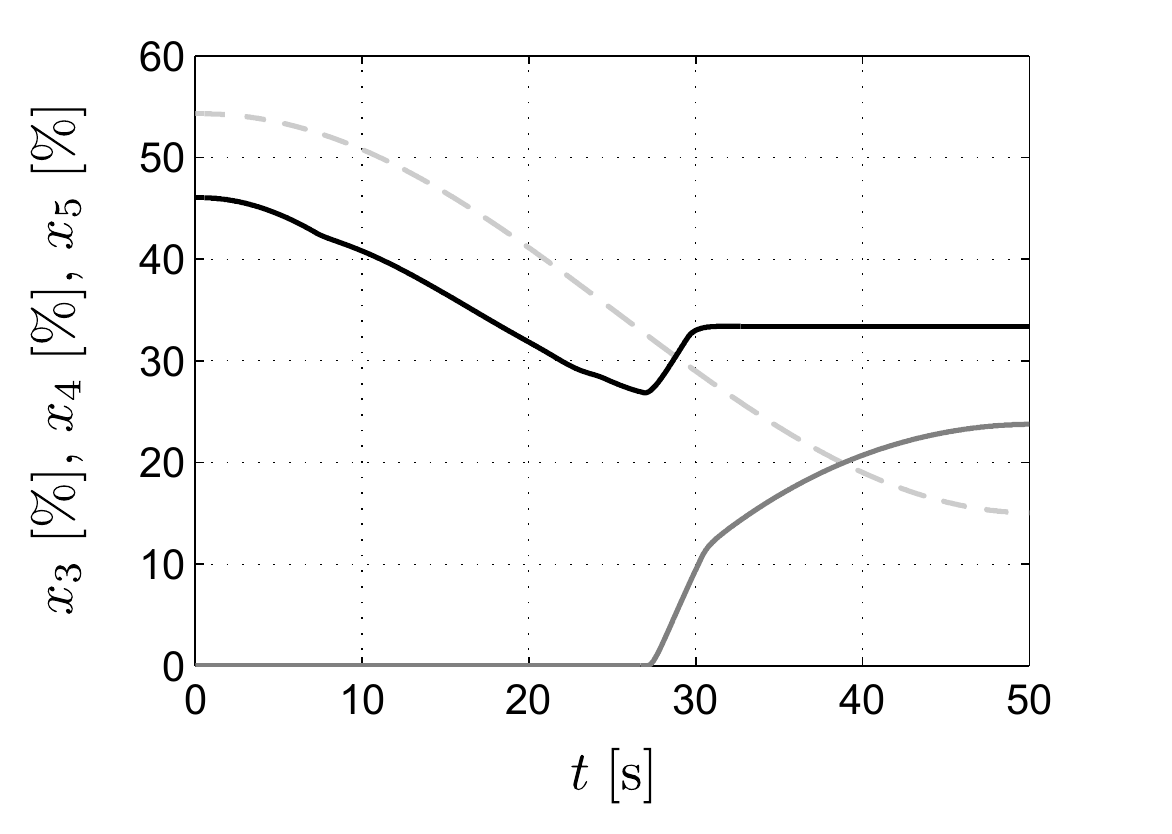}	\caption{Guide vane(GV) position $x_3$ (solid black),  process valve (PV) position $x_4$ (dashed light grey) and blow-off valve (BOV) position $x_5$ (solid dark grey) for scenario 2.} \label{fig:SimActuatorsS2}
        \end{subfigure}\hfill
    \end{subfigure}\hfill

	\caption{Simulation results of the pressure ratio $x_2$ and actuator positions $x_3$, $x_4$, $x_5$.}
	\label{fig:SimPiAndActuators}
\end{figure}

\begin{figure}[!h]
    \begin{subfigure}[t]{1\textwidth}
        \begin{subfigure}{0.485\textwidth}
        	\includegraphics[width=\linewidth]{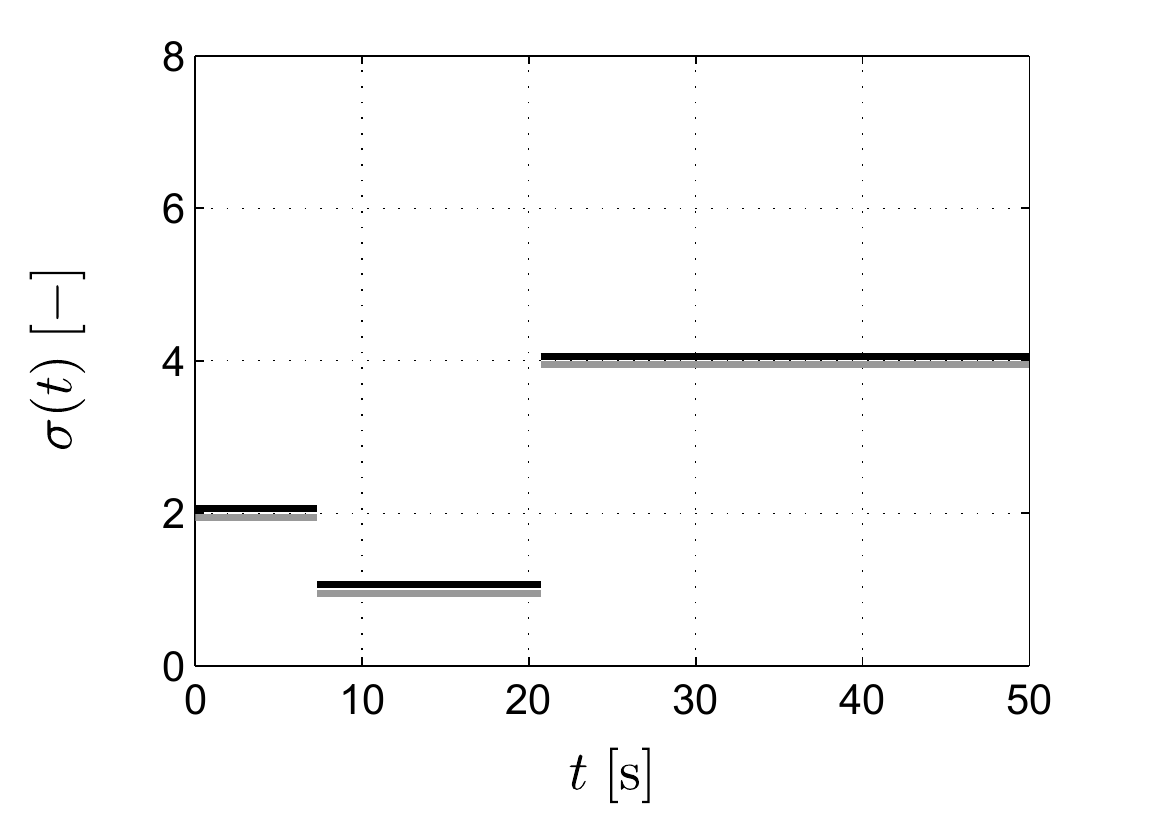}	\caption{Switching signal for guide vane controllers (solid black) and BOV controllers (grey) for scenario 1.}
            \label{fig:SimSwitchSigS1}
        \end{subfigure}\hfill
        \begin{subfigure}{0.485\textwidth}
            \includegraphics[width=\linewidth]{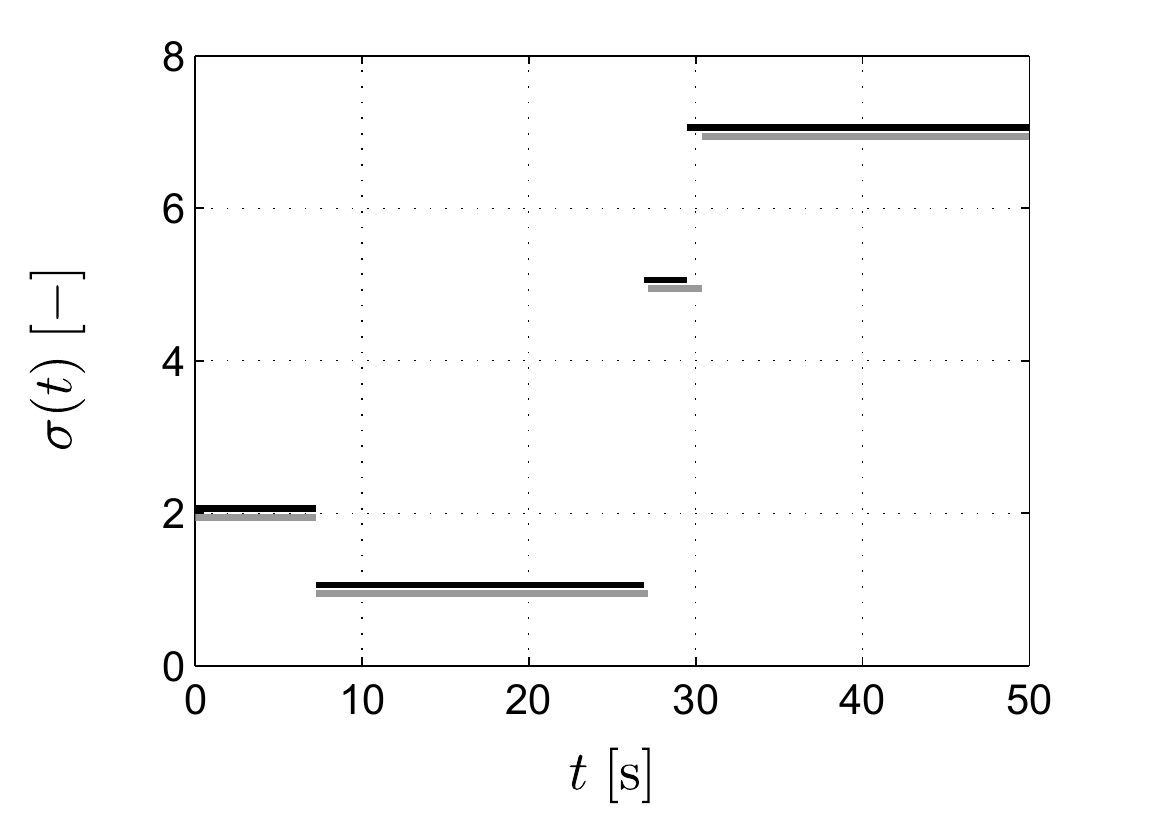}	\caption{Switching signal for guide vane controllers (solid black) and BOV controllers (grey) for scenario 2.}
            \label{fig:SimSwitchSigS2}
        \end{subfigure}\hfill
    \end{subfigure}\hfill

    \begin{subfigure}[t]{1\textwidth}
        \begin{subfigure}{0.485\textwidth}
        	\includegraphics[width=\linewidth]{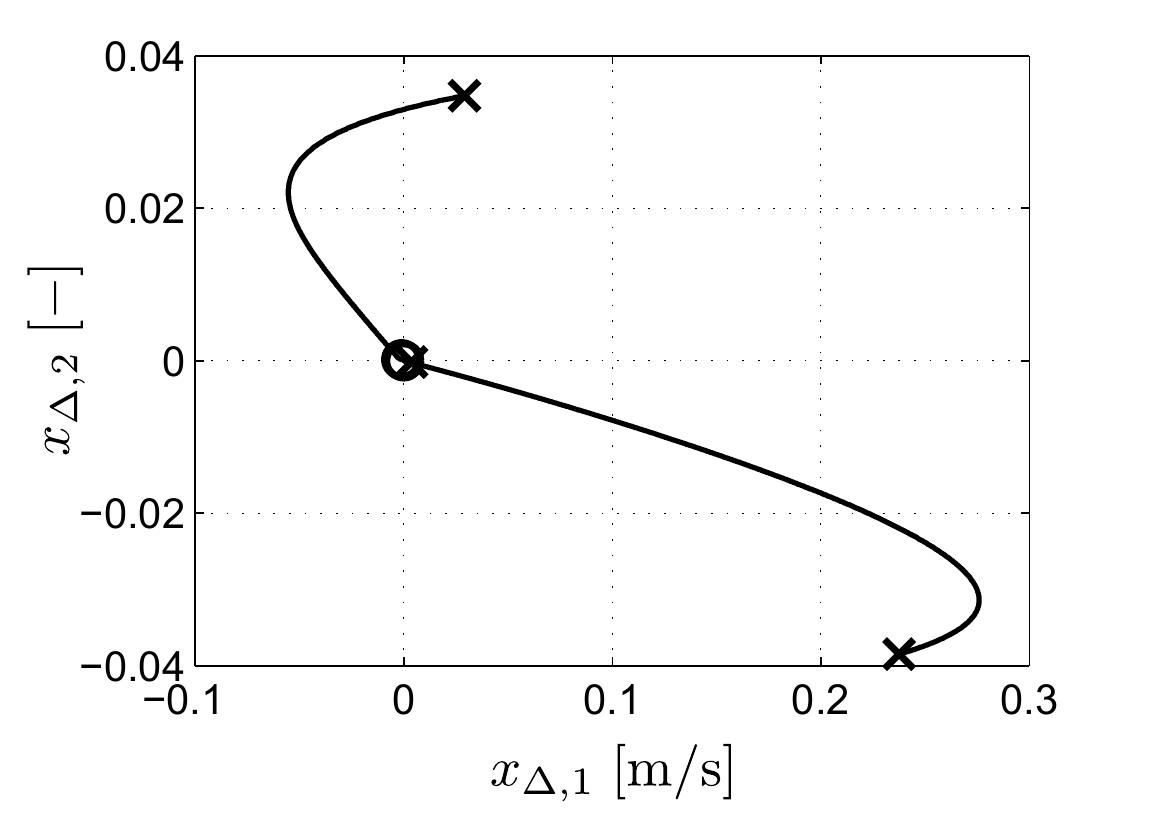}	\caption{Compressor working point motion displayed in $x_{\Delta,1}$-$x_{\Delta,2}$-plane (solid black) for scenario 1.}
            \label{fig:SimWPMotionES1}
        \end{subfigure}\hfill
        \begin{subfigure}{0.485\textwidth}
            \includegraphics[width=\linewidth]{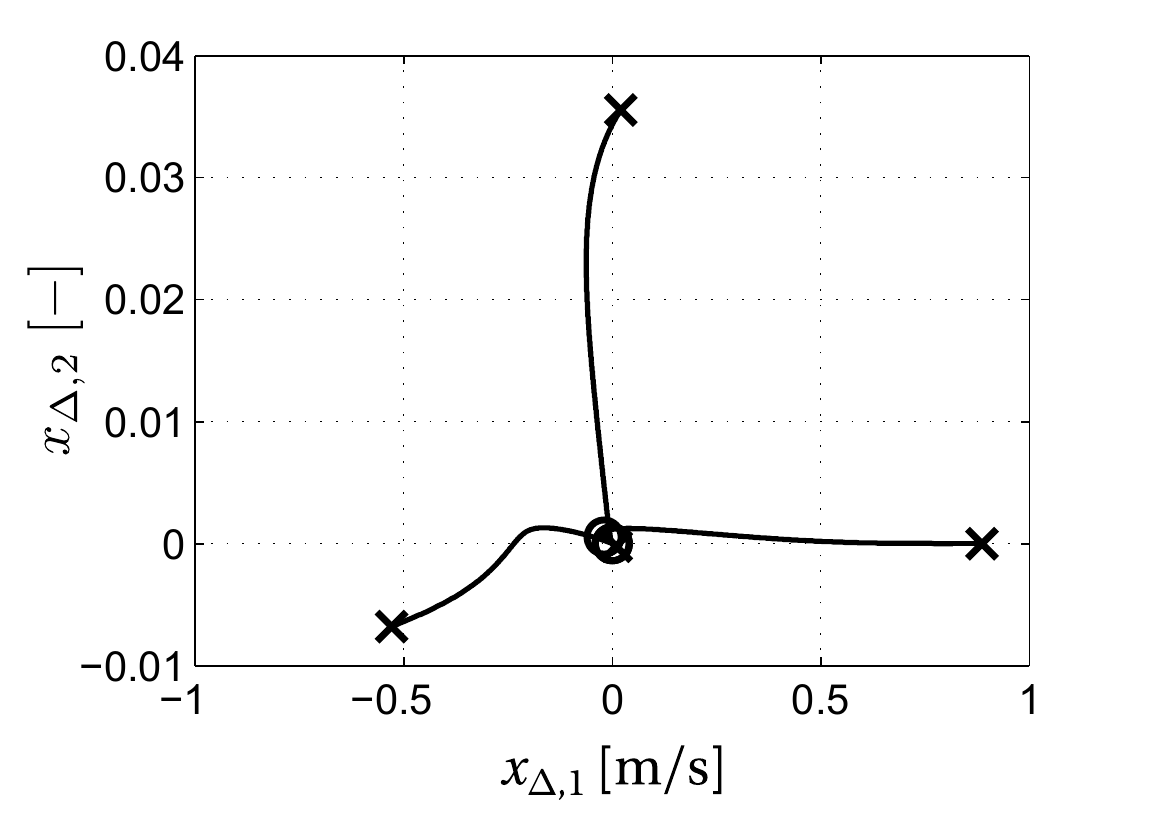}	\caption{Compressor working point motion displayed in $x_{\Delta,1}$-$x_{\Delta,2}$-plane (solid black) for scenario 2.}
            \label{fig:SimWPMotionES2}
        \end{subfigure}\hfill
    \end{subfigure}\hfill

	\caption{Simulation results of the switching signal $\sigma$ and compressor working point motion displayed in error coordinates.}
	\label{fig:SimSwitchSigAndWPMotionE}
\end{figure}

In this section we show the results of two simulation scenarios of the combined compressor system \eqref{eq:AutoSysWithControlExok} under MIN/MAX override control as presented in the last chapter. Thus, working principles of the control as well as successful output regulations are demonstrated. The combined compressor system \eqref{eq:AutoSysWithControlExok} uses the compressor model  \eqref{eq:ComprModelFinal} as plant description. The parameters of the model are given in Table \ref{tab:parameter}, see Appendix A.
In both scenarios the main control objective is to regulate the discharge pressure while simultaneously rejecting process disturbance introduced by actions of the process valve (PV). Thus, the main controller trajectory is designed to increase the pressure ratio from $x_2 = 1.7$ to $x_2 = 2.0$.  Scenario 1 is depicted in Figures~\ref{fig:SimWPMotionS1},~\ref{fig:SimPiS1},~\ref{fig:SimActuatorsS1},~\ref{fig:SimSwitchSigS1} and~\ref{fig:SimWPMotionES1}, while for scenario 2, see   Figures~\ref{fig:SimWPMotionS2},~\ref{fig:SimPiS2},~\ref{fig:SimActuatorsS2},~\ref{fig:SimSwitchSigS2} and~\ref{fig:SimWPMotionES2}.
The main difference between scenario 1 and scenario 2 is the severity of the process disturbance. 
While in scenario 1 the discharge pressure limit is set to $x_{2,max} = 1.9$, in scenario 2 it is set to $x_{2,max} = 2.05$, i.e. the discharge pressure limitation is virtually switched off in scenario 2. Although in both scenarios a process disturbance is introduced by the PV; see Figures~\ref{fig:SimActuatorsS1} and~\ref{fig:SimActuatorsS2}, the disturbance in scenario 2 is designed to be larger to drive the compressor towards the surge limit. Thus, the blow-off valve (BOV) becomes active.

According to the setup described above, the compressor will reach the maximum pressure limit in the scenario 1 (see Figure~\ref{fig:SimWPMotionS1}) while reaching the surge control line (SCL) in scenario 2 (see Figure~\ref{fig:SimWPMotionS2}). Simulation results show that the compressor system first follows the ramp-like main trajectory as shown in Figures~\ref{fig:SimPiS1} and~\ref{fig:SimPiS2} before switching to the corresponding limiting controller.
The switching can be  seen in Figures~\ref{fig:SimSwitchSigS1} and~\ref{fig:SimSwitchSigS2}. There are two switching events in  scenario 1. The first event corresponds to switching between the constant and sinusoidal main trajectory that have been introduced in Section~\ref{subsec:trajectory} and  \ref{subsec:appl} to approximate a ramp-like setpoint motion; i.e. $C_2 \rightarrowtail C_1$. The second event corresponds to $C_1 \rightarrowtail C_4$ which is the override by the maximum pressure discharge limiter. During the remaining simulation the system is limited to the maximum value of $x_{2,max} = 1.9$ while rejecting the process disturbance introduced by the PV. As illustrated by Figure~\ref{fig:SimActuatorsS1} the disturbance is rejected by adjusting the guide vane (GV) position only. Hence, the BOV is kept close and the compressor system is subjected to the controllers designed for the SISO domain only. As the BOV does not become active in scenario 1, the BOV position $x_5(t)=0$ for the entire simulation shown in Figure~\ref{fig:SimActuatorsS1}.

Keeping the BOV closed corresponds to a very energy efficient way of maximum pressure limitation. In common practice, the BOV is used for maximum pressure limitation due to its very fast response time. However, the BOV becoming active implies venting of already compressed gas and therefore a rather energy inefficient operation. Furthermore, Figure~\ref{fig:SimWPMotionES1} shows that the system is stable as defined in Section~\ref{ssec:stability}. Both switching events (indicated by black circles) take place in the close neighbourhood of $\tilde{x}=0$. Hence, the closed-loop system satisfies the regulation condition~\eqref{cond:ControlDesignRegulation} in this simulation.

In scenario 2 there is one more switching event; see Figure~\ref{fig:SimSwitchSigS2}. Similar to scenario 1 the simulation starts by following the ramp-like main trajectory motion; see Figure~\ref{fig:SimWPMotionS2}. As the PV is throttled much more, the system now reaches the SCL. Hence, there are overrides $C_1 \rightarrowtail C_5$ and $C_5 \rightarrowtail C_6$ by the anti-surge sub-controllers. These controllers now stabilize the system on the SCL while rejecting the process disturbance introduced by the PV. As the ASC is acting on the BOV, the system is now subjected to the sub-controllers designed for the MIMO domain. Finally, it can be seen from Figure~\ref{fig:SimWPMotionES2} that the closed-loop system is again practical stable and satisfies the regulation condition~\eqref{cond:ControlDesignRegulation} as desired during the controller design procedure.

\clearpage

\section{Conclusions} \label{Sec:Conclusions}

We considered an application-oriented nonlinear control of centrifugal compressors. The control is based on a modification of a model of Gravdahl and Egeland \cite{Gravdahl:1999b}. The modification 
accounts for characteristic features of industrial compressors and models the flow through the impeller and the diffuser separately. Thus, the total effective passage length depends on the pressure ratio  and facilitates for higher differential pressures. Furthermore, the modified model describes a typical scenario in process control of industrial compressors with control signals simultaneously and independently adjusting the positions of the guide vane and the blow-off valve while rejecting disturbances introduced by the process valve.

In the paper, we combined the framework of nonlinear output regulation via the internal model principle with MIN/MAX-override control. The overall design meets the requirements of industrial applications of centrifugal compressors, as the control is capable of tracking constant and time-dependent trajectories, as well as generating state constraints. A setpoint generator for time-varying setpoints has been designed using the proposed method of MIN/MAX-override control. This allows approximating a ramp-like setpoint motion  by  suitable overrides of the system responses generated from the corresponding marginally stable exo-systems. In industrial compressor operations we typically find surge avoidance rather than active suppression, which corresponds to the state constraints in the controller design, while tracking trajectories matches adjustment to variable process demands. Override control implies a switching scheme and we have given a proof of practical stability for the overall system. The analytical and simulation results show that the override control can be applied to discharge pressure control, anti-surge control and
maximum discharge pressure limitation.

Future work can go in several directions. Particularly relevant from an  application point of view is the implementation of the presented MIN/MAX-override concept in a real centrifugal air compressor system. In this context, further aspects should be investigated. These could include robustness against model and parameter uncertainties as well as using output feedback based on measurement. In this case, an essential aspect is knowledge about disturbance inputs acting on the system. Since the proposed concept also takes these disturbances into account, it must be clarified how these disturbances can be measured or modelled in a given application.

Further work could be done in extending the compressor model to incorporate additional fluid dynamical phenomena. This may include the following aspects: (1) considering changes in rotational speed, (2) accounting for different input and output conditions which would allow modeling of multi-stage compressors, and (3) improving the model of the compressor map, e.g. by using pseudo-invariant characteristics or applying  efficiency maps which would also account for  the compressor behaviour under variable gas conditions.

Finally, another further research direction is to augment control concepts as described in this paper by applications of computational intelligence methods for compressors \cite{hafaifa:2010,Schulze:2012,spindler:2021,wu:2012}. Thus, the amount of  information extracted from the compressor operation and available for analysis and subsequent utilization would be enhanced considerably which is a promising way for an advanced control performance.

\clearpage
\section*{Appendix}
\appendix
\section{Centrifugal compressor modelling}

\begin{figure}[!h]
\begin{center}
\includegraphics[width=8cm]{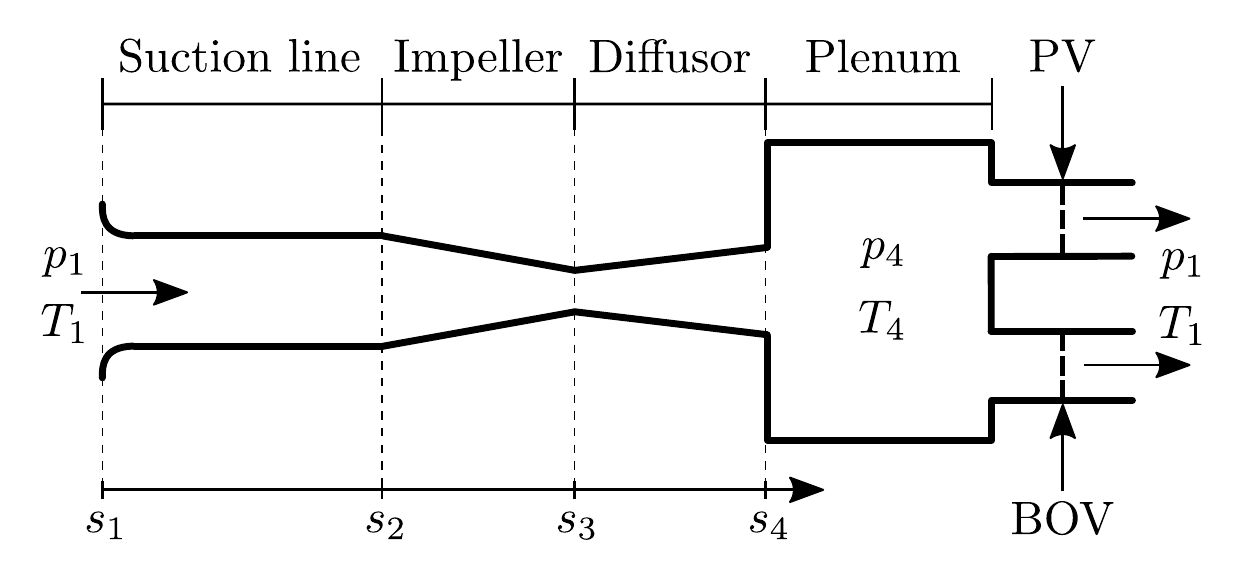}
\caption{Coordinate system for the centrifugal compressor model. Index $1$ refers to ambient conditions.}
\label{figSys}
\end{center}
\end{figure}

We propose a model which uses relevant physical dimensional parameters rather than dimensionless parameters, as they can be easily interpreted by practitioners. The modeling is done by presenting component models for each part of the system: suction line, impeller, diffuser, plenum, process valve (PV) and blow-off valve (BOV); see Figure~\ref{figSys} and Table~\ref{tab:Notation} for details. Comprising all component models will result in the complete model which will be considered for control design.

%Figure~(\ref{figMap}) shows an example of a compressor map depending on impeller's rotational speed. The dynamic behavior is described by ODEs or PDEs, depending on the used model.
%The printed column width is 8.4 cm. Size the figures accordingly.
%\begin{figure}[h!]
%\begin{center}
%\includegraphics[width=8cm]{KennfeldMitGrenzen.pdf}
%\caption{Stable region of a compressor map (pressure ratio $\pi$ vs. inlet speed $c$) depending on rotational speed $\omega$. The compressor operation is limited by surge limit to the left, choke limit to the right and maximum speed to the top.}
%\label{figMap}
%\end{center}
%\end{figure}

%The printed column width is 8.4 cm. Size the figures accordingly

\renewcommand{\arraystretch}{1.1}
\begin{table}[!h]
\centering
\scalebox{0.8}{
\begin{tabular}{cl}
\hline\hline
 Symbol & Description \\ \hline
 $s_i$ & Axial coordinate \\
 $c_i$ & Velocity of fluid at $s_i$ \\
 $p_i$ & Static pressure of fluid at $s_i$ \\
 $T_i$ & Absolute temperature of fluid at $s_i$ \\
 $\rho_i$ & Density of fluid at $s_i$ \\
 $A_i$ & Cross section of component at $s_i$ \\
 $L_{ij}$ & Covered distance by fluid particle from $s_i$ to $s_j$ \\
 $\Delta h_{ij}$ & Specific enthalpy difference from $s_i$ to $s_j$ \\
 $\Delta u_{ij}$ & Specific inner energy difference from $s_i$ to $s_j$ \\ 
 $\dot{m}_C$ & Compressor mass flow \\ 
 $\dot{m}_{PV}$ & Process valve mass flow \\
 $\dot{m}_{BOV}$ & Blow-off valve mass flow \\ \hline\hline \\
\end{tabular} }
\caption{Notation for the compressor model.}
\label{tab:Notation}
\end{table}
\renewcommand{\arraystretch}{1.5}

\begin{figure}[!h]
	\begin{subfigure}[t]{1\textwidth}
		\begin{subfigure}{0.485\textwidth}
			\includegraphics[width=\linewidth]{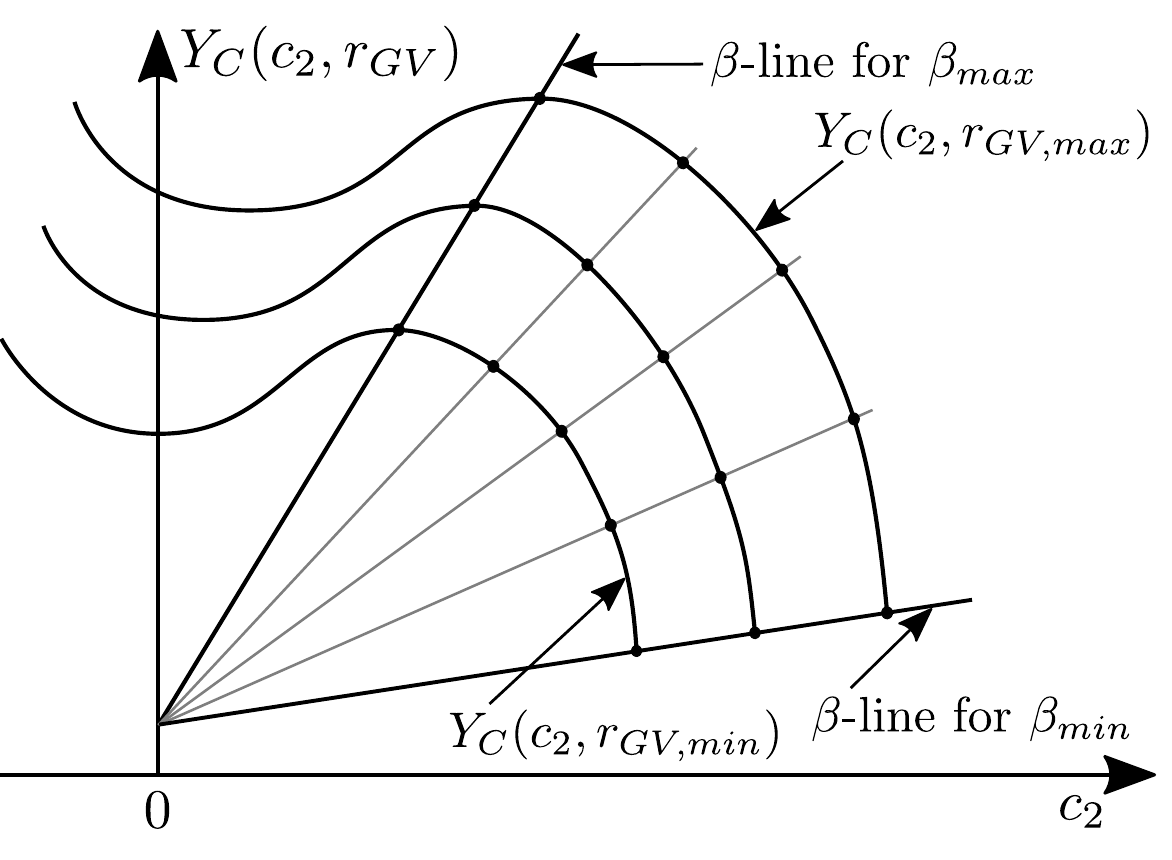}	\caption{Interpolation of the stable compressor map $Y_C(c_2, r_{GV})$ using some guide vane positions $r_{GV}$.} \label{ComprMapInterpolation}
		\end{subfigure}
		\hspace*{\fill}
		\begin{subfigure}{0.485\textwidth}
			\includegraphics[width=\linewidth]{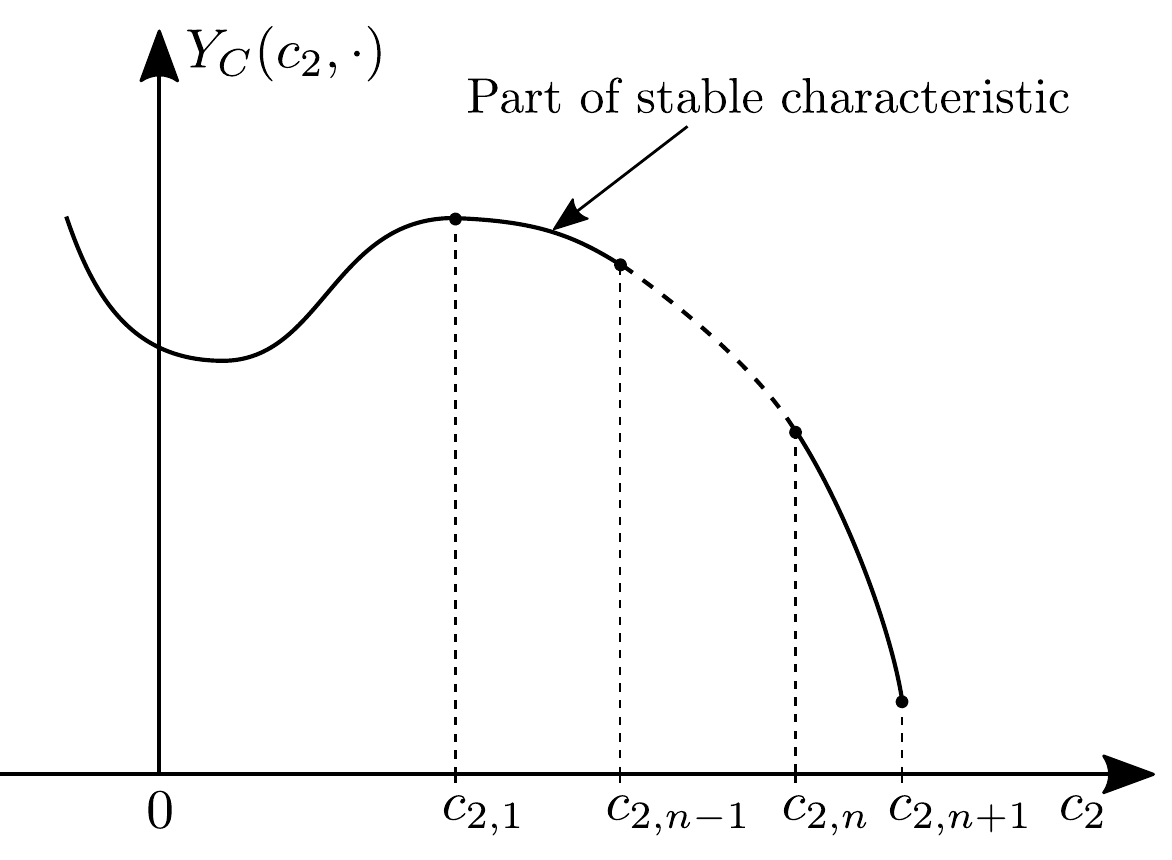}	\caption{Modeling the stable part of the compressor characteristic at constant guide vane positions $r_{GV}$.} \label{ComprCharacteristicPoly}
		\end{subfigure}
  \end{subfigure}\hfill
	\caption{Static compressor map modeling using $\beta$-line interpolation.}
	\label{fig:ComprCharacteristic}
\end{figure}

The suction line is described as one-dimensional, transient and threadlike incompressible flow without friction and heat transfer across piping boundaries. Therefore, we have $c_1 = c_2$. Now, integrating the corresponding transient Bernoulli-equation
\begin{equation}
	\frac{1}{\rho} \int_{p_1}^{p_2}dp =-\int_{s_1}^{s_2}\frac{\partial c}{\partial t}ds
	\label{eq:SuctionModel}
\end{equation}
gives:	
\begin{equation}
	\frac{p_2-p_1}{\rho_1} = -L_{12} \frac{\partial c_2}{\partial t}
	\label{eq:SuctionModelFinal}
\end{equation}
with $L_{12}$ being the effective length of the suction line. The centrifugal impeller does work on the system. Despite the complex flow field, the component model is based on a one-dimensional, steady-state and threadlike compressible flow that is derived from the first law of thermodynamics for open and stationary flow systems. The specific enthalpy difference $\Delta h_{23}$ is given by:
\begin{equation} \label{ImpellerModelSS}
	\Delta h_{23} = \Delta h_I + \frac{c^2_3 - c^2_2}{2}
\end{equation}
where $\Delta h_I$ is the specific enthalpy difference caused by the impeller. To account for the transient response of the compressor, we apply a deviation function that is proportional to the rate of change of the fluid's speed:
\begin{equation}
	\Delta h_{23} = \Delta h_I - \frac{c_3^2-c_2^2}{2} - \int_{s_2}^{s_3} \frac{\partial c}{\partial t}ds.
	\label{eq:ImpellerModelDyn}
\end{equation}
Comparable approaches have been proposed and successfully used by Moore \cite{Moore:1984}, Moore \& Greitzer \cite{Moore:1986b} and Spakovszky \cite{Spakovszky:2001}. Using enthalpy's definition and assuming a linear change of the cross sectional area-density ratio $A(s)/\rho(s)$ along $s$ as proposed by Spakovszky \cite{Spakovszky:2001}, the component model reads:
\begin{equation}
	\frac{p_3}{\rho_3} - \frac{p_2}{\rho_2} = \Delta h_I - \Delta u_{23} - \frac{c_3^2-c_2^2}{2} - L^*_{23}\frac{dc_2}{dt}ds
	\label{eq:ImpellerModel}
\end{equation}
where the effective impeller passage length is given by $L^*_{23}\mathrel{\mathop:}=L_{23}\ln(r_I)/(r_I-1)$ with $r_I\mathrel{\mathop:}=\rho_3 A_3/\rho_2 A_2$. The term $\Delta h_I - \Delta u_{23} - (c_3^2-c_2^2)/2$ can be interpreted as the steady state compressor map $Y_C$ of the centrifugal compressor leading to:
\begin{equation}
	\frac{p_3}{\rho_3} - \frac{p_2}{\rho_2} = Y_C - L^*_{23}\frac{dc_2}{dt}ds.
	\label{eq:ImpellerModelFinal}
\end{equation}
The diffuser collects all gas streams originating from impeller blade channels and transports them to the compressor outlet. The motion is  described as one-dimensional, transient and threadlike incompressible flow:
\begin{equation} \label{DiffuserFlow}
	\frac{1}{\rho}\int^{s_4}_{s_3}dp = - \frac{c^2_4 - c^2_3}{2} - \int^{s_4}_{s_3}\frac{dc}{dt}ds.
\end{equation}
No work is done on the flow and no heat is exchanged. The diffuser converts dynamic pressure into static pressure by increasing the cross section. Hence, the fluid velocity is decreased. We assume a linear change of the cross section along $s$. After performing the integration in~(\ref{DiffuserFlow}) this leads to the component model of the diffuser:
\begin{equation} \label{eq:DiffuserModelFinal}
	\frac{p_4-p_3}{\rho_4} = - \frac{c^2_4 - c^2_3}{2} - L^*_{34}\frac{dc_2}{dt}
\end{equation}
where the effective impeller passage length is given by $L^*_{34}\mathrel{\mathop:}=L_{34}\ln(r_D)/(r_I r_D-r_I)$ with $r_D\mathrel{\mathop:}=A_4/A_3$. The summation of Equations~\eqref{eq:SuctionModelFinal},~\eqref{eq:ImpellerModelFinal} and~\eqref{eq:DiffuserModelFinal} now gives the momentum equation from $s_1$ to $s_4$:
\begin{equation} \label{eq:CentrifugalCompModel}
	\frac{p_4}{\rho_4} - \frac{p_1}{\rho_1} = Y_C - L\frac{dc_2}{dt}
\end{equation}
with total effective passage length
\begin{equation} \label{eq:TotalEffectivePassageLength}
	L\mathrel{\mathop:}=L_{12} + L^*_{23}+L^*_{34}.
\end{equation}
The common procedure is to replace $Y_C$ in Equation \eqref{eq:CentrifugalCompModel} by the measured static compressor map describing the pressure increase depending on the compressor inlet velocity, i.e. $Y_C=Y_C(c_2)$. As we are considering guide vane controlled centrifugal compressors, the guide vane position $r_{GV}$ will be integrated, which is leading to the measured compressor map $Y_C(c_2, r_{GV})$. We are using an isoline approach for modeling the compressor map. The first type of isolines corresponds to the characteristic curve along the constant guide vane position $r_{GV}$. The second type of isolines corresponds to the so-called $\beta$-lines. The basic idea for the introduction of the $\beta$-lines is to create a bijective coordinate system \cite{Kurzke:1996}. This allows unique identification of the compressor state by the guide vane position $r_{GV}$ and the $\beta$-value. Furthermore, if the $\beta$-lines are chosen appropriately, each of these isolines can be assigned a physical interpretation, see Figure~\ref{ComprMapInterpolation}. For example, the $\beta_{max}$ line marks the transition between aerodynamically stable and aerodynamically unstable operation, i.e. the surge limit. The $\beta_{min}$ line marks the choke limit. Along the $\beta$-lines, the compressor exhibits similar aerodynamic behavior. The stable characteristic for constant guide vane position $r_{GV}$ will be modeled by cubic polynomials with continuous first derivative along $c_2$. Hence, the stable characteristic is $C^1$ for constant $r_{GV}$:
\begin{equation} \label{eq:MapModelComplete}
   Y_C(c_2,r_{GV}) =
   \begin{cases}
        a_{3,1} c_2^3 + a_{2,1} c_2^2 + a_{1,1} c_2 + a_{0,1} & c_{2,1} \le c_2 < c_{2,2},\\
        a_{3,2} c_2^3 + a_{2,2} c_2^2 + a_{1,2} c_2 + a_{0,2} & c_{2,2} \le c_2 < c_{2,3},\\
        \qquad\qquad\qquad\vdots & \qquad\quad\vdots \\
        a_{3,n} c_2^3 + a_{2,n} c_2^2 + a_{1,n} c_2 + a_{0,n} & c_{2,n} \le c_2 < c_{2,n+1}.
   \end{cases}
\end{equation}
where the index $n$ is the number of polynomials used to model the stable characteristic for a given $r_{GV}$, see Figure~\ref{ComprCharacteristicPoly}. The coefficients $a_{i,j}$ and interval boundaries $c_{i,j}$ are specific for a given $r_{GV}$, i.e. $a_{i,j}(r_{GV})$ and $c_{i,j}(r_{GV})$. The index $n=1$ is reserved for the surge limit of the compressor for a given $r_{GV}$.
Now, the dynamic behavior of the guide vane adjustment is considered. The guide vane motion is assumed to have first order low pass behavior. Hence, the dynamics of the guide vane position $r_{GV}$ is modeled by:
\begin{equation} \label{eq:DynModelGV}
	\frac{d r_{GV}}{d t} = \frac{1}{\tau_{GV}}\left[u_{GV} - r_{GV}\right].
\end{equation}
with the motion's time constant $\tau_{GV}$ and the control input $u_{GV}$.
\begin{figure}[!ht]
    \begin{center}
        \includegraphics[width=0.485\linewidth]{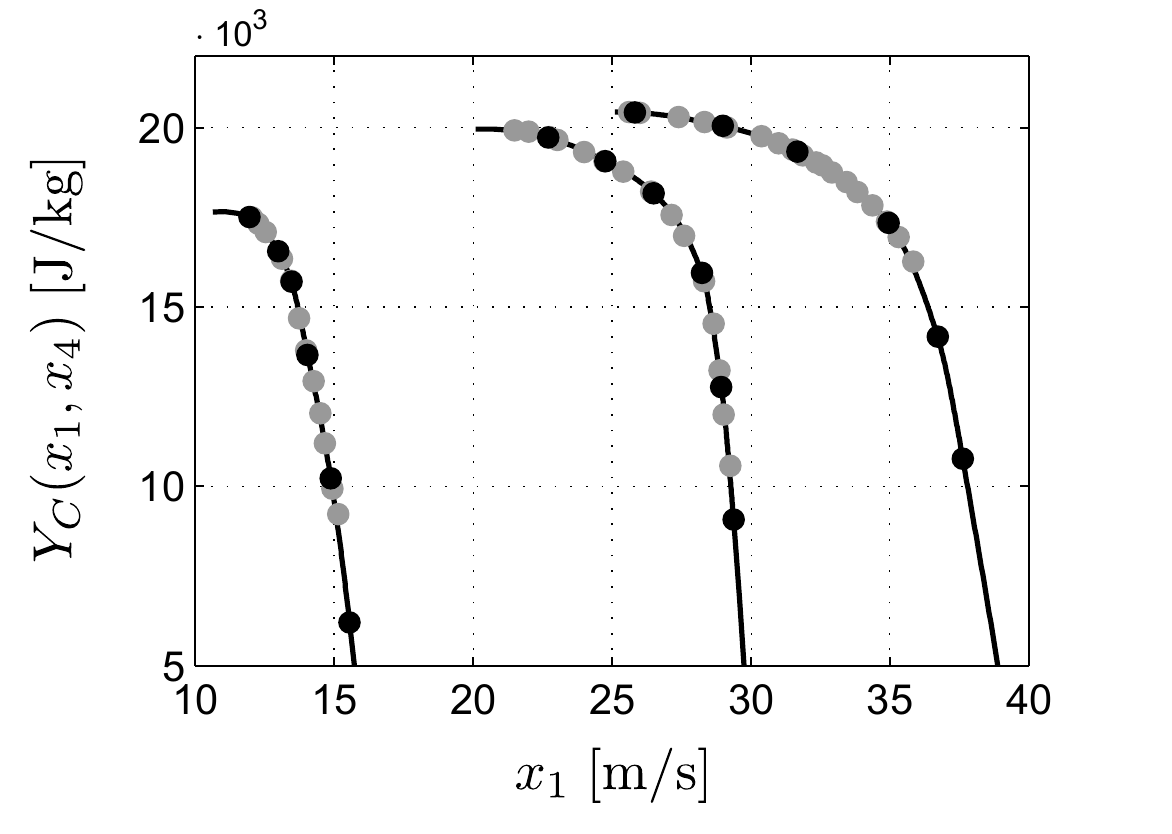}	
        \caption{Example map for an air compressor with outlet guide vane. The stable characteristics (black lines) are fitted to measured static compressor working points (gray dots). The transition points for polynomial interpolation are shown by gray dots.}
        \label{ComprMapOGV}
    \end{center}
\end{figure}
For simulations the measured compressor map from a real centrifugal air compressor with outlet guide vanes (OGV) has been used for parameterization. An example of the resulting map of a centrifugal air compressor can be seen in Figure~\ref{ComprMapOGV}. 
%Furthermore,  examination of $L$ reveals that its value only varies insignificantly compared to a suitably chosen mean value $\bar{L}$. The maximum deviation does not exceed $2.0\%$ within the inspected range of $p_4$. Hence, we will use the mean passage length $\bar{L}$ in Equation~(\ref{ImpulseEquation}). 

The modeling is completed by deriving models of the plenum and the valves. Using the conservation of mass and assuming an isentropic process, the plenum is described by:
\begin{equation} \label{PlenumModel}
	\frac{dp_4}{dt} = \frac{\kappa p_4}{V\rho_4}\left(\dot{m}_C - \dot{m}_{PV} - \dot{m}_{BOV}\right)
\end{equation}
where $V$ is the plenum volume and $\kappa$ is the isentropic exponent. The flows through the PV and the BOV are usually modeled using a quadratic approach, see for example~\cite{Gravdahl:1999b}. However, its validity refers to subcritical and incompressible valve flows only. Since centrifugal compressors can reach high pressure ratios, the energy equation for ideal, frictionless and compressible gases can be used to derive Bendemann's equation. This approach models the isentropic valve flow from a pressure vessel through a nozzle \cite{Sigloch:2005}:
\begin{subequations} \label{ThrottleModel}
	\begin{align}
        \dot{m}_{PV} &= A_{PV}\Psi \sqrt{\frac{2\kappa}{\kappa-1}\rho_4 p_4}, \\
		\dot{m}_{BOV} &= A_{BOV}\Psi \sqrt{\frac{2\kappa}{\kappa-1}\rho_4 p_4},
	\end{align}
\end{subequations}
with the cross section $A_{PV}$ and $A_{BOV}$ for process valve and blow-off valve as well as the Nusselt outlet function $\Psi$:
\begin{equation} \label{eq:NusseltOutflow}
	\Psi = \sqrt{\left(\frac{p_5}{p_4}\right)^{\frac{2}{\kappa}} - \;\left(\frac{p_5}{p_4}\right)^{\frac{\kappa+1}{\kappa}}}.
\end{equation}
Adjusting the PV and the BOV positions means changing the cross-sectional areas $A_{PV}$ and $A_{BOV}$.~The relative stroke $h_{PV} := A_{PV} /A_{PV,max}$ for the PV and the relative stroke $h_{BOV} := A_{BOV} /A_{BOV,max}$ for the BOV are introduced with the maximum cross-sectional areas $A_{PV,max}$ and $A_{BOV,max}$. However, when the valve positions $r_{PV}$ and $r_{BOV}$ are changed, the relative strokes $h_{PV}$ and $h_{BOV}$ are adjusted depending on specific valve designs and the resulting characteristics $Y_{PV}$ and $Y_{BOV}$. Thus, $h_{PV}=Y_{PV}(r_{PV})$ and $h_{BOV}=Y_{BOV}(r_{BOV})$ with $Y_{PV}(0)=A_{PV,min}$ and $Y_{PV}(1)=A_{PV,max}$ for the PV as well as $Y_{BOV}(0)=A_{BOV,min}$ and $Y_{BOV}(1)=A_{BOV,max}$ for the BOV. Now, considering the outlet function \eqref{eq:NusseltOutflow} and the valve characteristics $Y_{PV}$ as well as $Y_{BOV}$, the plenum model \eqref{PlenumModel} can be rewritten as:
\begin{equation} \label{eq:PlenumModelRe}
	\frac{d p_4}{d t} = \frac{\kappa A_2 p_4 \rho_2}{V \rho_4} \left[c_2 - \left[k_{PV}Y_{PV}(r_{PV})+k_{BOV}Y_{BOV}(r_{BOV})\right]\Psi\sqrt{\frac{2\kappa}{\kappa-1}\frac{\rho_4 p_4}{\rho_2^2}}\right].
\end{equation}
with:
\begin{subequations} \label{eq:ValveFlowParameters}
	\begin{align}
		k_{PV}(x_2) &= 
		\begin{cases}
			{\scriptstyle K_{PV}}\dfrac{{\scriptstyle A_{PV,\text{max}}}}{{\scriptstyle A_2}} & {\scriptstyle x_2 \le} \left(\dfrac{{\scriptstyle 2}}{{\scriptstyle \kappa+1}}\right)^{-\frac{{\scriptstyle 1}}{{\scriptstyle r_k}}},\\
			{\scriptstyle K_{PV}}\dfrac{{\scriptstyle A_{PV,\text{max}}}}{{\scriptstyle \sqrt{2}A_2}}\left(\dfrac{{\scriptstyle 2x_2^{r_k}}}{{\scriptstyle \kappa+1}}\right)^{\frac{{\scriptstyle \kappa+1}}{{\scriptstyle 2\kappa-2}}}\left(\dfrac{{\scriptstyle \kappa-1}}{{\scriptstyle x_2^{r_k}-1}}\right)^{\frac{{\scriptstyle 1}}{{\scriptstyle 2}}} & {\scriptstyle x_2 >} \left(\dfrac{{\scriptstyle 2}}{{\scriptstyle \kappa+1}}\right)^{-\frac{{\scriptstyle 1}}{{\scriptstyle r_k}}},
		\end{cases}\\
		k_{BOV}(x_2) &= 
		\begin{cases}
			{\scriptstyle K_{BOV}}\dfrac{{\scriptstyle A_{BOV,\text{max}}}}{{\scriptstyle A_2}} & {\scriptstyle x_2 \le} \left(\dfrac{{\scriptstyle 2}}{{\scriptstyle \kappa+1}}\right)^{-\frac{{\scriptstyle 1}}{{\scriptstyle r_k}}},\\
			{\scriptstyle K_{BOV}}\dfrac{{\scriptstyle A_{BOV,\text{max}}}}{{\scriptstyle \sqrt{2}A_2}}\left(\dfrac{{\scriptstyle 2x_2^{r_k}}}{{\scriptstyle \kappa+1}}\right)^{\frac{{\scriptstyle \kappa+1}}{{\scriptstyle 2\kappa-2}}}\left(\dfrac{{\scriptstyle \kappa-1}}{{\scriptstyle x_2^{r_k}-1}}\right)^{\frac{1}{2}} &  {\scriptstyle x_2 >} \left(\dfrac{{\scriptstyle 2}}{{\scriptstyle \kappa+1} }\right)^{-\frac{{\scriptstyle 1}}{{\scriptstyle r_k}}},
		\end{cases}
	\end{align}
\end{subequations}
where
$k_{PV}=K_{PV}A_{PV,max}/A_2$ and $k_{BOV}=K_{BOV}A_{BOV,max}/A_2$ are for subcritical valve flow. $K_{PV}$ and $K_{BOV}$ are correction factors valid for a speciﬁc valve. In the simulations we use an equal percentage characteristics for the PV with $Y_{PV}$ given by:
\begin{equation} \label{eq:PVCharactericts}
    Y_{PV}(r_{PV}) = k_{PV,0}\left(\dfrac{1}{k_{PV,0}}\right)^{r_{PV}}
\end{equation}
with $k_{PV,0}$ being the dimensionless $K\hspace{-0.7mm}v_0$ design parameter that can be extracted from the corresponding data sheet. For the BOV we use a linear characteristic.

The process valve position $r_{PV}$ and the blow-off valve position $r_{BOV}$ is controlled by the inputs $u_{PV}$ and $u_{BOV}$. The valve motion is assumed to have first order low pass behavior. Hence, the dynamics of the PV position $r_{GV}$ and the dynamics of the BOV position $r_{BOV}$ are modeled by:
 
	\begin{align}
        \frac{d r_{PV}}{d t} &= \frac{1}{\tau_{PV}}\left[u_{PV} - r_{PV}\right], \label{eq:DynValveModels} \\
		\frac{d r_{BOV}}{d t} &= \frac{1}{\tau_{BOV}}\left[u_{BOV} - r_{BOV}\right]. \label{eq:DynValveModels1} 
	\end{align}
Equations \eqref{eq:CentrifugalCompModel}, \eqref{eq:DynModelGV}, \eqref{eq:PlenumModelRe}, \eqref{eq:DynValveModels}  and \eqref{eq:DynValveModels1} constitute the complete compressor system model. Now, these equations are rearranged to contain only the density ratio $\rho_4/\rho_2$, the pressure ratios $p_4/p_2$ and $p_5 /p_4$ as well as the pressure-density ratios $p_1/\rho_1$ and $p_2/\rho_2$. These ratios are subsequently replaced by known quantities. If the compressor works on ambient upstream and downstream conditions, then $p_5 = p_1$. If friction is neglected across the suction piping then $p_1 = p_2$. This reduces the number of required ratios. The model now depends only on the density ratio $\rho_4/\rho_1$, the pressure-density ratio $p_1/\rho_1$ and the pressure ratio $\Pi\mathrel{\mathop:}=p_4/p_1$. The pressure-density ratio $p_1/\rho_1$ can be described by the equation of state (EOS) for ideal gases $p_1/\rho_1 = R_S T_1$, where $R_S$ is the speciﬁc gas constant and $T_1$ is the ambient temperature. Both are assumed to be known. Now, we use the pressure ratio $\Pi=p_4/p_1$ describing the relation between plenum pressure and ambient pressure as a further state variable and replace the density ratio and the pressure-density ratio $\rho_4/\rho_1$. For this we assume an isentropic process for perfect gas during impeller and diffuser passage:
\begin{equation} \label{eq:IsentropicSuct2Dis}
	\frac{\rho_4}{\rho_1} = \left(\frac{p_4}{p_1}\right)^{\frac{1}{\kappa}}= \Pi^{\frac{1}{\kappa}}.
\end{equation}
 Thus, the total effective passage length $L$ in Equation \eqref{eq:TotalEffectivePassageLength} becomes:
\begin{equation} 
L(x_2) = L_{12} + \frac{L_{23}}{\Pi^{\frac{1}{\kappa}}\frac{A_3}{A_2}-1}\ln\left(\Pi^{\frac{1}{\kappa}}\frac{A_3}{A_2}\right) + \frac{L_{34}}{\Pi^{\frac{1}{\kappa}}\frac{A_3}{A_2}{\left(\frac{A_4}{A_3}-1\right)}}\ln\left(\frac{A_4}{A_3}\right),
\end{equation}
which completes the modeling.
The model can now be reformulated in its final version using the state vector $x=(c_2,\Pi,r_{GV},r_{PV},r_{BOV})^T$; see Equation~\eqref{eq:ComprModelFinal}. In the numerical simulations presented in Section \ref{subsec:sim}, the model parameters given in Table \ref{tab:parameter} are used.

\renewcommand{\arraystretch}{1.1}
\begin{table}[!h]
\centering
\scalebox{0.8}{
\begin{tabular}{lll}
\hline\hline
 Description & Symbol &   Value  \\ \hline
 Specific gas constant [$J K^{-1} mol^{-1}$] & $R_S$ & 287.0 \\ 
 Ambient temperature [K] & $T_1$ & 295.4 \\ 
 Isentropic exponent [-] & $\kappa$ & 1.4 \\ 
 Effective length of suction line [m] & $L_{12}$ & 13.0\\
 Effective length of impeller [m] & $L_{23}$ & 2.5 \\
 Effective length of diffuser [m] & $L_{34}$ & 1.0 \\
 Cross section of impeller [$m^2$] & $A_2$ & 0.44 \\
 Cross section of diffuser inlet [$m^2$] & $A_3$ & 0.11 \\
 Cross section of diffuser outlet [$m^2$] & $A_4$ & 0.22 \\  
 Plenum volume [$m^3$] & $V$ & 32.0 \\
 Process valve correction factor [-] & $K_{PV}$ & 1.0 \\
 Process valve maximum cross-sectional area [$m^2$] & $A_{PV,max}$ & 0.196 \\
 Process valve minimum flow parameter [-] & $K\hspace{-0.7mm}v_0$ & 0.03 \\
 Blow-off valve correction factor [-] & $K_{BOV}$ & 1.0 \\
 Blow-off valve maximum cross-sectional area [$m^2$] & $A_{BOV,max}$ & 0.196 \\
 Time constant of guide vane motion [s]& $\tau_{GV}$ & 0.50 \\
 Time constant of process valve motion [s]& $\tau_{PV}$ & 0.35 \\
  Time constant of blow-off valve motion [s]& $\tau_{BOV}$ & 0.35 \\
 \hline\hline \\
\end{tabular} }
\caption{Parameter of the compressor model.}
\label{tab:parameter}
\end{table}
\renewcommand{\arraystretch}{1.5}

\end{document}